\documentclass[11pt]{amsart}
\usepackage{stmaryrd}
\usepackage{cases}
\usepackage{mathrsfs}
\usepackage{txfonts}
\usepackage{bbm}
\usepackage{amsfonts}
\usepackage{dsfont}
\usepackage{graphicx}
\usepackage{amstext}
\usepackage{amssymb}
\usepackage{latexsym}
\usepackage{subfigure}
\usepackage{tikz}
\usepackage{ulem}
\usepackage{CJK}
\usepackage{float}
\usepackage[linktocpage=true]{hyperref}
\usepackage{listings}

\def \nn{\nonumber}
\def \I{{\mathrm I}}
\def \II{{\mathrm {II}}}
\def \III{{\mathrm {III}}}

\def \R{\mathbb R}
\def \Z{\mathbb Z}
\newfont{\bbb}{cmb11 at 8pt}
\def \beq{\begin{equation}}
\def \eeq{\end{equation}}
\def \beqs{\begin{equation*}}
\def \eeqs{\end{equation*}}
\def \bea{\begin{eqnarray}}
\def \eea{\end{eqnarray}}
\def \beas{\begin{eqnarray*}}
\def \eeas{\end{eqnarray*}}
\def \brk{\begin{remark}}
\def \erk{\end{remark}}

\newenvironment{prof}[1][\textbf{Proof.}]{\noindent\textit{#1}\quad } {\hfill$\Box$\vspace{0.7mm}}
\newenvironment{prof1}[1][\textbf{Proof of Theorem \ref{Holder}.}]{\noindent\textit{#1}\quad } {\hfill$\Box$\vspace{0.7mm}}
\newenvironment{prof2}[1][\textbf{Proof of Lemma \ref{lemma.LDT.5}.} ]{\noindent\textit{#1}\quad } {\hfill$\Box$\vspace{0.7mm}}
\makeatletter

\newcommand{\Rmnum}[1]{\expandafter\@slowromancap\romannumeral #1@}
\makeatother
\newcommand{\Keywords}[1]{\par\noindent {\small{ Keywords\/}: #1}}

\newtheorem{definition}{Definition}[section]
\newtheorem{theorem}{Theorem}[section]
\newtheorem{lemma}{Lemma}[section]
\newtheorem{corollary}{Corollary}[section]
\newtheorem{remark}{Remark}[section]

\begin{document}

\title[H\"{o}lder continuity of Lyapunov exponent]{H\"{o}lder continuity of Lyapunov exponent for a family of smooth Schr\"{o}dinger cocycles}\thanks{This work was supported by NSFC of China (Grants 11471155, 11771205).}
\author{Jinhao Liang, Yiqian Wang, Jiangong You}
\address{Department of Mathematics, Southeast University, Nanjing 210096, China}\email{jinhaohao999@163.com}
\address{Department of Mathematics, Nanjing University, Nanjing 210093, China}\email{yiqianw@nju.edu.cn}
\address{Chern Institute of Mathematics and LPMC, Nankai University, Tianjin 300071, China} \email{jyou@nju.edu.cn}

\begin{abstract}
We prove the H\"older continuity of the Lyapunov exponent for quasi-periodic Schr\"odinger cocycles with a $C^2$ cos-type potential and any fixed Liouvillean frequency, provided the coupling constant is sufficiently large . Moreover, the H\"older exponent is independent of the frequency and the coupling constant.

\

\Keywords{Schr\"{o}dinger cocycle; Lyapunov exponent; Liouvillean frequency; H\"older continuity.}
\end{abstract}

\maketitle

\section{Introduction}

Consider the one-dimensional quasi-periodic Schr\"odinger operator on $l^2(\mathbb Z)$
\beq\label{operator}
(H_{\alpha,\lambda,v,x}u)_n=u_{n+1}+u_{n-1}+\lambda v(T^nx)u_n,
\eeq
$$T:\ x \longmapsto x+\alpha.$$
We call $v$, $\lambda$, $x$, $\alpha$ potential ,coupling constant, phase and frequency respectively.
%Here $v\in C^2\left(\R/\Z, \mathbb R\right)$ is the potential, $\lambda\in \mathbb R$ the large coupling constant, $x\in \R/\Z$ the phase, and $\alpha \in \R/\Z$ the irrational frequency.
Let $A^{(E-\lambda v)}\in C^2\left(\R/\Z, \text{SL}(2,\mathbb R)\right)$ be given by
\beq\label{cocycle}
A^{(E-\lambda v)}(x)=\left(
  \begin{array}{cc}
    E-\lambda v(x) & -1 \\
    1              & 0  \\
  \end{array}
\right).
\eeq
 Then
$$
(x,w) \mapsto \left(x+\alpha, A^{(E-\lambda v)}(x)w\right)
$$
defines a family of dynamic systems on $\R/\Z\times \mathbb R^2$, denoted for simplicity by $(\alpha,A^{(E-\lambda v)})$.
We call them the Schr\"{o}dinger cocycles associated with the Schr\"odinger operator (\ref{operator}).
The $n$th iteration of the cocycle is denoted by
$$(\alpha,A^{(E-\lambda v)})^n=(n\alpha, A_n^{(E-\lambda v)}),$$
where
$$
A_n^{(E-\lambda v)}(x)=\left\{
                         \begin{array}{ll}
                           A^{(E-\lambda v)}(x+(n-1)\alpha)\cdots A^{(E-\lambda v)}(x), & \hbox{ $n\geq1$;} \\
                           I_2, & \hbox{ $n=0$;} \\
                           {\left(A_{-n}^{(E-\lambda v)}(x+n\alpha)\right)}^{-1}, & \hbox{ $n\leq-1$.}
                         \end{array}
                       \right.
$$
The Lyapunov exponent $L(E,\lambda)$ of the cocycle is defined as
$$
L(E,\lambda)= \lim_{n\rightarrow\infty}\frac 1n \int_{\R/\Z} \log\|A_n^{(E-\lambda v)}(x)\|dx,
$$
where $\|\cdot\|$ denotes the matrix norm on $\text{SL}(2,\mathbb R)$. The limit always exists by Kingman's subadditive Ergodic Theorem. It is known that, for irrational $\alpha$,
$$
L(E,\lambda)=\lim_{n\rightarrow\infty}\frac 1n \log\|A_n^{(E-\lambda v)}(x)\|,
\quad\ \ \text{a.e.}\ \  x\in \R/\Z.
$$

Besides the Lyapunov exponent (short for LE), another important subject is the integrated density of states (short for IDS), which is given by
$$N(E)=\lim_{n\rightarrow \infty}\frac1n\int_{\mathbb R/\Z}\#\{(-\infty, E)\cap \Sigma(H_n(x))\}dx,$$
where $\Sigma(H_n(x))$ denotes the set of eigenvalues of the operator $H(x)$ restricted to $[1,n]$ with Dirichlet boundary conditions and $\#$ denotes the cardinality of a set.
It is well known that the IDS always exists and is locally constant outside the spectrum \cite{AS}.
Moreover, the IDS and the LE are related via the famous Thouless' Formula
$$L(E)=\int_{\mathbb R} \log|E-E'|dN(E').$$
By the Hilbert transform and the theory of singular integral operators, the H\"older continuity passes from $L(E)$ to $N(E)$ and vice versa (see \cite{GS1} for more details).

\subsection{The main result}
In this paper, we focus on a class of potentials, which has the following properties:
\begin{itemize}
  \item Assume $\frac {dv}{dx}=0$ at exactly two points, one is minimal and the other maximal, which are denoted by $z_1$ and $z_2$;
  \item These two extreme points are non-degenerate, i.e. $\frac{d^2v}{dx^2}(z_j)\neq 0$ for $j=1,2$.
\end{itemize}
We call this kind of $v$ is $C^2 $ $\cos$-type function.

For any irrational $\alpha $, let $\{{p_n}/{q_n}\}_{n\in \mathbb N}$ be the fraction approximation of $\alpha$. Define
\beqs
\beta(\alpha) :=\limsup_{n\rightarrow\infty} \frac{\log q_{n+1}}{q_n}.
\eeqs
We call $\alpha $ Liouvillean  if $0\leq\beta<\infty$ and we call $\alpha $ strong Liouvillean if $\beta=\infty$.
We also call $\alpha$ Diophantine if there exist $\tau\ge 1$ and $c>0$ such that \beqs
q_{n+1}\leq c q_n^\tau,\quad n\geq 1.
\eeqs
It is easy to see that every Diophantine number is  Liouvillean by our terminology, since we have $\beta(\alpha)=0$ for any Diophantine $\alpha$.
It is also clear that the set of Diophantine numbers has a full Lebesgue measure on $\mathbb R$.

Our main result is as follows.
\begin{theorem}\label{Holder}
Let $v$ be a $C^2$ $\cos$-type function and $\alpha$ a  Liouvillean number. Then there exists a $\lambda_1=\lambda_1\left(\alpha, v\right) $  such that for $\lambda\geq \lambda_1$, the LE and the IDS of the cocycles (\ref{cocycle}) are H\"older continuous in $E$, i.e. for $E,E'\in [E_1,E_2]$,
\beqs
\left|L(E)-L(E')\right|+\left|N(E)-N(E')\right|<C\left|E-E'\right|^\sigma,
\eeqs
where the constant $C$ depends on the interval $[E_1,E_2]$, and the H\"older exponent $\sigma$ is an absolute constant  depends on nothing and is bigger than $10^{-8}$.
\end{theorem}

\brk
For $C^2$ $cos$-type potentials considered in Theorem \ref{Holder}, Anderson Localization has been established by Sinai \cite{S} and Fr\"ohlich-Spencer-Wittwer \cite{FSW}. In \cite{FSW}, they also assumed that the potentials are even functions.
Uniform positivity of LE was also proved by different methods, see Bjerkl\"ov \cite{Bjer}, Wang-Zhang \cite{WZ} and Liang-Kung \cite{LK}.
\erk

\begin{remark} When the frequency is Diophantine, Wang-Zhang \cite{WZ} obtained a weaker regularity ($\log$-H\"older continuity) of the LE for quasi-periodic Schr\"odinger cocycles with $C^2$ $\cos$-type potentials. Theorem \ref{Holder} improves their result and provides the first positive example on the H\"older continuity of LE on the energy $E$ for finite smooth potentials and fixed Liouvillean frequencies.

\end{remark}

\begin{remark}\label{sigma}
The H\"older exponent $\sigma$ in Theorem \ref{Holder} is independent of $\lambda$. In contrast, the H\"older exponent obtained in \cite{GS1,YZ} will tend to zero as $\lambda$ goes to infinity. Furthermore, the lower bound of $\sigma$ can be improved by a more careful estimate.
However, we do not attempt to obtain the optimal value here.
\end{remark}

%\brk
%Damanik-Lukic-Goldstein \cite{DLG} showed that the Schr\"odinger operator has homogenous spectrum if the IDS is H\"older continuity and the spectral gaps decay exponentially. It is hopeful using our methods to obtain the asymptotic property of the gaps for $C^2$ $\cos$-type potentials.
%\erk

\subsection{Related results and some remarks}
In recent years, much work has been devoted to the regularity properties of the LE and the IDS for the one-dimensional discrete quasi-periodic Schr\"odinger operators (cocycles).
Goldstein-Schlag \cite{GS1} developed the technique from Bourgain-Goldstein \cite{BG}, and obtained a sharp version of large deviation theorem (short for LDT) for general real analytic potentials with strong Diophantine frequencies.
They developed the Avalanche Principle to prove H\"older continuity of LE in the regime of positive LE. However, the H\"older exponent they obtained tends to zero as $\lambda$ goes to infinity.
This contradicts the intuition that the H\"older exponent should get better when increasing $\lambda$.
Later You-Zhang \cite{YZ} improved the method in \cite{BJ} and \cite{B1}, and extended \cite{GS1}'s result to all Diophantine frequencies and some Liouvillean frequencies. Nevertheless, the H\"older exponent in \cite{YZ} also goes to zero when $\lambda\rightarrow \infty$.
Very recently, Han-Zhang \cite{HZ} showed the H\"older exponent is in fact independent of %the Lyapunov exponent (and hence of $\lambda$)
$\lambda$, and extended \cite{YZ}'s result to all Liouvillean frequencies.

The results of regularity of the LE on trigonometric polynomial potential are much better. Bourgain \cite{B2} obtained $(\frac12-\varepsilon)$-H\"older regularity of LE for almost Mathieu operator with large coupling.
He noted that the result seems optimal, that is, one should not expect to have global H\"older regularity better than $\frac12$.
Goldstein-Schlag \cite{GS2} showed that IDS is $(\frac1{2k}-\varepsilon)$-H\"older continuous for small analytic perturbations of trigonometric polynomials of degree $k$ and positive LEs and strong Diophantine frequencies.
%They investigated the distribution of the Dirichlet eigenvalues of $H(x,\alpha)\varphi=E\varphi$ on a finite interval $[1,N]$ and obtained the lower bound for the separation of the eigenvalues, which contributes to the H\"older regularity.
Recently, Ge-You \cite{GY}
improved the H\"older exponent to $\frac1{2k}$ for trigonometric polynomials of degree $k$.
We would like to emphasize that all the results stated above require positivity of LE or sufficiently large coupling constants.

It seems that the regularity of IDS or LE is better for small coupling.
Amor \cite{Am} used the sophisticated KAM iteration in \cite{E} to prove the $\frac12$-H\"older continuity for quasiperiodic Schr\"odinger operators with Diophantine multi-frequencies, where the smallness of the coupling depends on the frequency.
{ Avila-Jitomirskaya \cite{AJ} got the $\frac12$-H\"older exponent of the IDS with small coupling nonperturbatively in subexponential regime (i.e. $\alpha$ satisfies $\beta(\alpha)=0$).
Their result was extended to Liouvillean frequencies in \cite{LY}.}
We remark that so far all the approaches heavily depend on the analyticity of the cocycles, which are not applicable to $C^2$ smooth cocycles as considered in this paper.

For rough potentials, only a few results are established for the continuity of LE.
Klein \cite{K} obtained the $\log$-H\"older
continuity of the LE for some Gevrey potentials and strong Diophantine frequencies. In \cite{Chan}, Chan used multi-scale analysis and obtained uniform positivity of LE for some $C^3$ potentials by excluding a positive measure of frequencies and by varying the potentials in some typical way.
Wang-Zhang \cite{WZ} obtained the uniform positivity and $\log$-H\"older continuity of LE provided that the frequency is any fixed Diophantine number, the coupling is large and the potential is a $C^2$ $\cos$-type function.
This is the first result on finite smooth potentials.
%More precisely, they developed a new method from \cite{Y,Z1}, which is in spirit to the techniques in \cite{BC}, to overcome the difficulty of resonance (see Section 2 for rigorous definition) and establish an induction scheme, from which a weak version of LDT follows.
%However, this kind of LDT is not enough to prove H\"older regularity.

For other related results, Avila-Krikorian \cite{AK} studied so-called monotonic cocycles which are a class of smooth or analytic cocycles non-homotopic to identity. They prove that the LE is smooth or even analytic respect to the parameter.
Jitomirskaya-Kachkovskiy \cite{JK} showed the IDS is Lipschitz continuous for the Schr\"odinger cocycles with strictly monotone Lipschitz potentials and all irrational frequencies.
In comparison, the global regularity of LE is considered to be not better than $\frac 12$-H\"older continuous for quasi-periodic Schr\"odinger cocycles.

The key for this paper is to improve the large deviation theorem in \cite{WZ}. In \cite{WZ}, the upper bound for the measure of the set of phases with a large deviation is sub-exponentially small. To guarantee the H\"older regularity, we need to prove that the upper bound is in fact exponentially small.

Some new technical difficulties appear when the frequency is extended to Liouvillean numbers.
A Liouvillean frequency may lead to a sudden loss of the finite Lyapunov exponent.
%It was solved in \cite{LK}, where a Liouvillean version of the induction theorem is established to obtain the estimates of monodromy matrices whenever the phase belongs to the critical intervals.
Moreover, it may cause some complicated cases, where we have to deal with the resonances in different scales.
%To avoid this, we will choose a subsequence of the fraction approximation of the frequency, which is inspired by the so-called "CD"-bridge technique in \cite{AFK}, to simplify the procedure.

We remark that the relation between $\beta$ and $\lambda$ in Theorem \ref{Holder} is roughly
\beq\label{lnlambdabeta}
\log\lambda\geq C(\beta+1),\eeq
where  $C>0$ is independent of $\beta$.
Thus if $\beta\rightarrow \infty$, we have to require $\lambda\rightarrow \infty$.
Hence our method is invalid for the case with strong Liouvillean frequencies.
However, we guess in general there is no H\"older continuity when $\beta=\infty$.
\

The rest of the paper is organized as follows.
In Section 2, we will first give some technical lemmas, which are proved in \cite{WZ}, to estimate the norms and the angles of the monodromy matrices.
And then we will review the induction theorem with Liouvillean frequency in \cite{LK} .
In Section 3 we will show the sharp large deviation theorem.
It is the key part of the paper.
We will ensure the exponential growth of the norm of $A^{(E-\lambda v)}_n(x)$ excluding a exponentially small set of $x$.
In Section 4 we will prove a technical lemma used in Section 3.
%Since our result concerns about Liouvillean frequency, it requires us to carefully analyze the trajectory of the base dynamic.
In Section 5 we will finish the proof of the H\"older continuity using the LDT and the Avalanche Principle.

\section{Preliminaries}
\subsection{Lemmas for product of matrices}
In this subsection we present several technical lemmas to compute the norms of the product of some matrices. All of them are proved in \cite{WZ}.

For the Schr\"odinger cocycle (\ref{operator}), we have $\|A^{(E-\lambda v)}(x)\|\geq 1$. The equality may hold for some $x$ and $E$, which will cause some technical difficulties.
The next lemma transforms the Schr\"odinger cocycle (\ref{operator}) to a cocycle $(\alpha,A(x))$ satisfing $\|A(x)\|>1$ for all $x\in \R/\Z$.

\begin{lemma}\label{lemma.transform}
Let $J$ be any compact interval. For $x\in\mathbb R/\mathbb Z$, and $t=E/\lambda\in J$, define the following cocycles map
\beqs
A(x)=\Lambda(x)\cdot R_{\phi(x,t)}:=\left(\begin{array}{cc} \lambda(x) & 0 \\ 0 & \lambda^{-1}(x) \\ \end{array} \right)\cdot  \left(\begin{array}{cc} \cos \phi(x,t) & -\sin \phi(x,t)  \\ \sin\phi(x,t) & \cos\phi(x,t) \\ \end{array} \right).
\eeqs
 Assume also that $\lambda(x)$ satisfies
$$ \lambda\leq \lambda(x)<C\lambda, \ \ \left|\frac{d^m\lambda(x)}{dx^m}\right|<C\lambda, \quad m=1,2.
$$
Then $A^{(E-\lambda v)}(x)$ is conjugate to $A(x)$, where $\cot \phi(x,t)\rightarrow t-v(x)$ in $C^2$-topology as $\lambda\rightarrow \infty$.
Therefore for $\lambda$ sufficiently large(only depends on $v$), to prove Theorem \ref{Holder}, it is enough to prove the corresponding results for $A(x)$ with $\cot \phi(x,t)= t-v(x)$.
\end{lemma}

Let $A\in \text{SL}(2,\mathbb R)$ satisfying $\|A\|>1$.
Define the map
$ s,u:\text {SL}(2,\mathbb R)\rightarrow \mathbb {RP}^1:=\mathbb R/ (\pi\mathbb Z)$,
such that $s(A)$ is the most contraction direction of $A$ and $u(A)=s(A^{-1})$.
In other words, let $\hat s (A)\in s(A)$ be the unit vector, we have $\|A\cdot\hat s(A)\|=\|A\|^{-1}$, and so is $u$.
We say $s(A)$ is the stable direction of $A$ and $u(A)$ is the unstable direction of $A$.
\iffalse For $\theta\in \mathbb R/\mathbb Z$, let
\beas R_{\theta}=\left(
                             \begin{array}{cc}
                               \cos(2\pi\theta) & -\sin(2\pi\theta) \\
                               \sin(2\pi\theta) & \cos(2\pi\theta) \\
                             \end{array}
                           \right)
                           \in \text{SO}(2,\mathbb R).
\eeas\fi
Then for $A\in \text{SL}(2,\mathbb R)$ with $\|A\|>1$, we have the polar decomposition
\bea\label{polardecomposition}
                        A=R_u\left(
                             \begin{array}{cc}
                               \|A\| & 0 \\
                               0 & \|A\|^{-1} \\
                             \end{array}
                           \right)R_{\frac \pi 2-s},
\eea
where $s,u\in [0,2\pi)$ are some suitable choices of angles corresponding to the directions $s(A), u(A) \in \mathbb {RP}^1$. For convenience, we also use $s(A)$ and $u(A)$ for $s$ and $u$.

Consider the product $E(x)=E_2(x)\cdot E_1(x)$. By (\ref{polardecomposition}), we have
\begin{small}
\beas
E=E_2\cdot E_1=R_{u(E_2)}\left(\begin{array}{cc}\|E_2\| & 0 \\ 0 & \|E_2\|^{-1} \\
                       \end{array}
                \right)R_{\frac \pi 2-s(E_2)+u(E_1)}
                \left(\begin{array}{cc}\|E_1\| & 0 \\ 0 & \|E_1\|^{-1} \\
                      \end{array}
                \right)R_{\frac \pi 2-s(E_1)}.
\eeas
\end{small}
We say $\left|s(E_2)-u(E_1)\right|$ is the angle between $E_2$ and $E_1$.
For $\|E_i\|\gg 1, \ i=1, 2$, it is not hard to see that the norm of $E_2E_1$ approximately  equals to $\left\|E_2\right\|\cdot\left\|E_1\right\|$ unless the angle tends to $0$.
We say it is nonresonant if
$$
\left|s(E_2)-u(E_1)\right|^{-1}\ll \min\left\{\|E_1\|,\|E_2\|\right\}.
$$
Otherwise, we say it is resonant.
The following lemma gives a basic idea to compute $s(E)$ and $u(E)$.

\begin{lemma}\label{lemma.basic}
Let
$$s(x)=s(E(x)),\quad u(x)=u(E(x)),\quad \theta(x)=s(E_2(x))-u(E_1(x)),$$
$$e_j(x)=\|E_j(x)\|,\quad s_j(x)=s_j(E_j(x)),\quad u_j(x)=u(E_j(x)), \quad j=1,2,$$
$$f_1=\frac12(e_1^2\tan\theta+e_1^2e_2^{-4}\cot\theta),\quad f_2=\frac12(e_2^2\tan\theta+e_2^2e_1^{-4}\cot\theta),$$
$$0<\kappa \leq 10^{-3},\quad 0< \eta \leq10^{-2}.$$
Then for each $m=0,1,2$ and each $x\in I$, we have the following.
\begin{itemize}
  \item If $e_2^\kappa\geq e_1\gg 1$, then we have
  \begin{small}$$
  c<\left|\frac{d^m(s-s_1)}{dx^m}/\frac{d^m{(\tan^{-1}(e_1^2\tan\theta)-\frac\pi2)}}{dx^m}\right|, \ \ \left|\frac{d^m(u-u_2)}{dx^m}/\frac{d^m{\cot^{-1}(\sqrt{f_2^2+1}+f_2)}}{dx^m}\right|<C.$$
  \end{small}
  \item If $e_1^\kappa\geq e_2\gg 1$, then we have
  \begin{small}$$c<\left|\frac{d^m(s-s_1)}{dx^m}/\frac{d^m{(\tan^{-1}(\sqrt{f_1^2+1}+f_1)-\frac\pi2)}}{dx^m}\right|, \ \
  \left|\frac{d^m(u-u_2)}{dx^m}/\frac{d^m{\cot^{-1}(e_2^2\tan\theta)}}{dx^m}\right|<C.$$
  \end{small}
  \item If $$e_2,\ e_1\gg 1,\quad \inf_{x\in I}|\theta|\geq \max_{x\in I}\{e_1^{-\eta},e_2^{-\eta}\},$$ then we have
      $$\|E\|\geq C e_1e_2|\theta|, $$ and
  \begin{small}$$c<\left|\frac{d^m(s-s_1)}{dx^m}/\frac{d^m{\tan^{-1}(e_1^2\tan\theta)}}{dx^m}\right|, \ \
   \left|\frac{d^m(u-u_2)}{dx^m}/\frac{d^m{\cot^{-1}(e_2^2\tan\theta)}}{dx^m}\right|<C.$$
  \end{small}
\end{itemize}
\end{lemma}
\begin{remark}\label{kappa}
In \cite{WZ}, $\kappa$ and $\eta$ need to be sufficiently small.
Straightforward computations show that $\kappa= 10^{-3}$ and $\eta = 10^{-2}$ are enough.
Moreover, the value $10^{-3}$ is not optimal and one can show that Lemma \ref{lemma.basic} still holds if $\kappa$ is somehow larger.
However, we do not want to obtain the optimal value here since it will be very complicated.
It is worth mentioning that the value of $\kappa$ is related to the constant $\delta$ in Theorem \ref{LDT}, and hence the H\"older exponent in Theorem \ref{Holder}.
%({\bf Reserve $\kappa$ in the proof of LDT and H\"older.})
\end{remark}
\begin{remark}
The first and the second claim in Lemma \ref{lemma.basic} are for the resonant case. As we see in the lemma, the stable direction of the product matrix differs greatly from the stable direction of the factor matrix. More precisely, if $\|E_2\|\gg\|E_1\|\gg 1$, then
$$s(E_2E_1)=s(E_1)+\arctan(\|E_1\|^2\tan(s(E_2)-u(E_1)))-\frac\pi2.$$
\end{remark}
\begin{remark}\label{2.3}
The third claim in Lemma \ref{lemma.basic} is for the nonresonant case.
One can repeat this argument to obtain the lower bound of the norm of the product of a sequence of SL$(2,\mathbb R)$ matrices.
In other words, if the angle $\theta_j$ between $E_{j+1}$ and $E_{j}$ is not close to $\frac\pi 2$ for $j=1,\cdots,n$, then the norm of their product has a lower bound
$$\|\prod_{j=1}^{n+1}E_{j}\|\geq C\cdot\prod_{j=1}^{n+1}\|E_{j}\|\cdot \prod_{j=1}^{n}|\theta_j|,$$
which is similar to the Avalanche Principle (see Lemma \ref{AP}).
Moreover, the stable direction of the product matrix approximately equals to the stable direction of the first factor matrix
$$\|s(\prod_{j=1}^{n+1}E_{j})-s(E_1)\|_{C^2}<\|E_1\|^{-(2-5\eta)}.$$
So is the unstable direction.
\end{remark}

\subsection{The induction theorem with Liouvillean frequency}
In \cite{WZ}, Wang and Zhang developed an induction scheme to overcome the difficulty of resonance and obtain the positivity of LE for $C^2$ $\cos$-type potentials and fixed Diophantine frequencies.
Liang-Kung \cite{LK} extended \cite{WZ}'s result to the Liouvillean frequencies.
They proved
\begin{theorem}\label{LZpositivity}\textbf{(\cite{LK})}
Consider the Schr\"odinger cocycle (\ref{cocycle}) with a $C^2$ $\cos$-type potential and a Liouvillean frequency. Then for $\varepsilon\in(0,10^{-3}]$, there exists a $\lambda_2=\lambda_2(\alpha,v ,\varepsilon)>0$ such that
$$L(E,\lambda)>(1-\varepsilon)\log\lambda$$
for all $(E,\lambda)\in \mathbb R \times [\lambda_2,\infty)$.
\end{theorem}

This theorem bases on the induction scheme as follows.
For $N$ large enough, we define the initial angle function $g_N$ as $g_N(x)=s(A(x))-u(A(x))=\arctan(t-v(x))$, where $A$ has the form in Lemma \ref{lemma.transform}, and the initial critical points $\left\{c_{N,1},c_{N,2}\right\}$ as $$\{c_{N,1},c_{N,2}\}= \{ y: |g_{N}(y)| = \min _{x\in \mathbb R / \mathbb Z}|g_{N}(x)| \}.$$
For $i\geq N$, we inductively define
\begin{itemize}
  \item The critical intervals $I_{i,1}$ and $I_{i,2}$ center at two critical points $\left\{c_{i,1},c_{i,2}\right\}$, with radii of $q_{i}^{-2}$ on $\R/\Z$.
  \item The return time $r_i^\pm(x):I_i\rightarrow \mathbb{Z}_+$ is the first return time after $10^{-3}\varepsilon q_{i}$, where $r_i^+(x)$ is the forward return and $r_i^-(x)$ backward. That is
      $$r_i^{\pm}(x)=\min\{j:T^{\pm j}x\in I_i,\ j\geq 10^{-3}\varepsilon q_{i}\},\quad x\in I_i.$$
      Let $r_i^\pm=\min_{x\in I_i} r_i^\pm(x)$ and $r_i=\min\{r_i^+,r_i^-\}$.
  \item The angle function $g_{i+1}: I_i\rightarrow \R/\Z$ is defined by $g_{i+1}(x)= s_{r_i^{+}}(x)-u_{r_i^{-}}(x) $, where $s_{r_i^{+}}(x)=s(A_{r_i^{+}}(x))$ and $u_{r_i^{-}}(x)=s(A_{-r_i^{-}}(x))$. Moreover, the critical points $\left\{c_{i+1,1},c_{i+1,2}\right\}$ are the minimal points of the angle function $g_{i+1}$ as the following
      \beas C^{(i+1)}=\left\{c_{i+1,1},c_{i+1,2}\right\} \subseteq \{ y: |g_{i+1}(y)| = \min _{x\in I_{i}}|g_{i+1}(x)| \}, \eeas
      with $$ |c_{i+1,j}-c_{i,j}| <C\lambda^{-\frac34r_{i-1}}, \ \ j=1,2. $$
\end{itemize}

\iffalse Recall (\ref{def.beta}).
We choose $N_0$ such that
$
\beta>\frac{1}{2}\sup_{n>N} q_n^{-1}\log q_{n+1}
$ for $N>N_0$.
For any $0<\varepsilon<1$, choose $N_1$ larger than $N_0$ such that
$
\sum_{n\geq N} q_n^{-1}<10^{-7}\varepsilon^2
$
for $N>N_1$. Fix $N>N_1$. Define
\beas
\log\lambda_n= \left\{
 \begin{array}{ll}
  \left(1-10^{-3}\varepsilon\right)\log \lambda, & \hbox{$n=N+1$;}\\
 \left(1-10^{5}\varepsilon^{-1}q_{n-1}^{-1}-10^{5}\varepsilon^{-1}q_{n}^{-1}\log q_{n}\right)\log \lambda_{n-1}, & \hbox {$n\geq N+2$.}
 \end{array}
\right.
\eeas
Straightforward calculation shows that $\lambda_n$ decreases to a $\lambda_\infty$ with $\lambda_{\infty}\geq \lambda^{1-3\cdot10^{-1}\varepsilon}$.
We choose $\lambda$ large enough such that
\beq
\label{thm.lambda0}
\lambda >\lambda_3(\alpha,v ,\varepsilon):=\max \left\{{C(v)}^{10\varepsilon^{-1}}q_N^{20\varepsilon^{-1}}, C(v)^{-100}q_N^{200}, e^{10^{9}\varepsilon^{-1}\beta}\right\}.
\eeq
\fi

\begin{theorem}\label{theorem.induction}\textbf{(\cite{LK})}
Let $v$,  $\varepsilon$, and $\alpha$  be as in Theorem \ref{LZpositivity}, and let $\lambda$ satisfy
\beq\label{thm.lambda0}
\lambda>\lambda_3(\alpha,v,\varepsilon):=\max \left\{C(v,\varepsilon,q_N), e^{10^{9}\varepsilon^{-1}\beta}\right\}.
\eeq
Then we have
$$\|A_{\pm r^\pm_i}(x)\| \geq \lambda^{(1-10^{-2}\varepsilon)r^\pm_i},\quad x\in I_{i}.$$
And the angle function $g_{i+1}$ satisfies the nondegenerate condition
\beq\label{nondegenerate}|g_{i+1}(x)|\geq C\|x-C^{(i+1)}\|^3, \quad x\in I_i \backslash B(C^{(i)},\lambda^{-r_i}),\eeq
where $\|x-C^{(i+1)}\|:=\min_{j=1,2}\|x-c_{i+1,j}\|$. Moreover, depending on the positions of the critical intervals (points), it is divided into three cases: Case $(i+1)_{\I}$, Case $(i+1)_{\II}$ and Case $(i+1)_{\III}$.
\begin{itemize}
\item In Case $(i+1)_\I$, it holds $$(I_{i,1}+j\alpha)\cap I_{i,2}=\varnothing,\quad j\in (-10^{-3}\varepsilon q_{i},10^{-3}\varepsilon q_{i})\cap \mathbb Z.$$
    This is a nonresonant case and we have that
    \beq\label{gi+1-gi} \|g_{i+1}-g_i\|_{C^2}\leq C\lambda^{-\frac32r_{i-1}}\quad \text{in} \ \ I_i.\eeq
\item In Case $(i+1)_\II$, the critical interval $I_i$ is a consecutive interval, i.e., $I_{i,1}\cap I_{i,2}\neq \varnothing$. This is also a nonresonant case and (\ref{gi+1-gi}) holds.
\item In Case $(i+1)_{\III}$, it holds
    $$(I_{i,1}+k\alpha)\cap I_{i,2}\neq\varnothing,\quad \exists \ k\in (-10^{-3}\varepsilon q_{i},10^{-3}\varepsilon q_{i})\cap \mathbb Z \backslash\{0\}.$$
    This is a resonant case and
    \beq\label{gi+1-gi3} g_{i+1}(x)=\arctan(l_k^2 \tan(g_{i}(x\pm k\alpha)))-\frac\pi2+g_{i}(x), \quad x\in I_{i},\eeq
    where $l_k=\|A_{k}(x)\|$, and "$\pm$" takes "$+$" for $x\in I_{i,1}$ and takes "$-$" for $x\in I_{i,2}$.
\end{itemize}
\end{theorem}

\begin{remark}\label{2.5}
For convenience, we assume $\varepsilon=10^{-3}$ in the rest of the paper. We will choose $\lambda_1=\lambda_3(\alpha,v,10^{-3})$ in Theorem \ref{Holder}, which implies (\ref{lnlambdabeta}).
\end{remark}

The proof of the induction theorem mainly follows \cite{WZ}.
However, the method in \cite{WZ} depends heavily on the arithmetic properties of Diophantine frequency.
When the frequency is Liouvillean, directly applying their method will cause a great loss of the norm of the $n$th iterated cocycle in the inductive step,
which might not be sufficient to get the positivity of the LE.
To reduce and control the loss, Liang-Kung gave delicate analysis on the base dynamic. From it they found that
the product of the angles is actually a factorial by the fact that the Liouvillean frequency can be well approximated by rational numbers.
Hence this product will be not too small, see Remark \ref{2.3}.
It was also found in \cite{LK} that each period  of resonance are followed by a much longer period of non-resonance, which can counteract the loss of finite Lyapunov exponent due to resonance.

%{\bf Give some ideas of the proof of this induction theorem.}

\section{The large derivation theorem}
In order to prove H\"older continuity of the Lyapunov exponent, Goldstein and Schlag \cite{GS1} established the so-called sharp large deviation theorem saying that
\beq\label{GSLDT}
\text{meas} \ \left\{ x\in \R/\Z: \left|\frac1n\log\|A_n(x)\|-L_n(E)\right|>\kappa\right\}<e^{-c\kappa n},
\eeq
where $L_n(E):=\frac1n\int_{\R/\Z}\log\left\|A_n(x)\right\|dx$,
provided the analyticity of the potentials, some Diophantine conditions of the frequencies and the positivity of the LE.
Combining (\ref{GSLDT}) with the Avalanche Principle (see Lemma \ref{AP}),
they obtained
the H\"older continuity of the LE.

For rough potentials, we intend to establish some kind of large deviation inequality.
More precisely, we have to ensure the exponential growth of the norm of $A^{(E-\lambda v)}_n(x)$ excluding a exponentially small set of $x$.
However, Wang-You \cite{WY1,WY2} gave some examples for general smooth Schr\"odinger cocycles that the norm of $A^{(E-\lambda v)}_n(x)$ may be not large for a large set of $x$ and the associated LE can be not continuous even for Diophantine frequencies.
Fortunately, by restricting the cocycle to the Schr\"odinger cocycle with $C^2$ $\cos$-type potential, Wang-Zhang \cite{WZ} proved the $\log$-H\"older continuity of the LE for any fixed Diophantine frequencies.
In fact, their proof depends on the following large deviation estimate
\beq\label{WZLDT}
\text{meas} \ \left\{ x\in \R/\Z: \left|\frac1n\log\|A_n(x)\|-L(E)\right| >\kappa\log\lambda \right\}<e^{-c\delta(\kappa) n^\sigma},
\eeq
for $\lambda>\lambda_0(\alpha,v,\kappa)$ and  $n>n_0(\alpha,v,\kappa)$, where $\kappa$, $\delta(\kappa)$ are positive numbers and  $\sigma$ satisfies $0<\sigma=\sigma(\alpha)<1$.
But this type of large deviation inequality is not enough to prove the H\"older continuity since the coefficient $\sigma$ in (\ref{WZLDT}) is strictly smaller than 1, see \cite {B1, GS1, YZ}.
Our main task in this paper is to replace $\sigma$ by 1 for all Liouvillean frequencies.
Now we are going to state our large deviation theorem.
\begin{theorem}\label{LDT}
Choose $N$ and $\lambda$ large enough (such that Theorem \ref{theorem.induction} holds). Let $i\geq i_0(\lambda,q_N)$.
Then the large deviation inequality holds as follows
\beq\label{LDT.LDT}
\text{meas} \ \left\{x\in\mathbb{R}/\mathbb{Z}:\frac 1 i \log\|A_i(x)\|<\frac{9}{10}\log \lambda\right\}<\lambda^{-\delta i},
\eeq
where the coefficient $\delta$ is independent of $\lambda$ and $\alpha$, and satisfies $\delta\geq10^{-6}$.
\end{theorem}

\begin{remark}
In fact, the coefficient $\delta$ is related to $\sigma$ in Theorem \ref{Holder} and $\kappa$ in Lemma \ref{lemma.basic}. For convenience, we will choose $\sigma=10^{-2}\delta=10^{-5}\kappa=10^{-8}$ in the following proof.
\end{remark}

%We first outline the main ideas of the proof of Theorem \ref{LDT}.
\begin{prof}To prove Theorem \ref{LDT}, the central problem is to obtain the lower bound of $\|A_n(x)\|$.
Note that by Theorem \ref{theorem.induction} whenever $x$ belongs to a critical interval $I_{j}$,  $\|A_{\pm r^\pm_{j}}(x)\|$ has a lower bound
$$\|A_{\pm r^\pm_{j}}(x)\|\geq \lambda^{(1-10^{-5})r^\pm_{j}},\quad x\in I_{j}.$$
We need to analyze the trajectory of the base dynamic and judge when $T^kx\in I_{j}$ for some suitable $j$ and $k$.
For this purpose, we divide the trajectory of the base dynamic in the following way.

\begin{figure}[H]\label{keylemma.figure}
\begin{center}
\begin{tikzpicture}[yscale=1]
\draw [densely dotted, thick] (-4,0) -- (-3,0);
\draw [thick] (-3,0) -- (3,0);
\draw [densely dotted, thick] (3,0) -- (4,0);
\draw plot [mark=*, mark size= 0.2ex] (-4,0);
\node[align=left, below] at (-4,0){$-l_3$};
\draw plot [mark=*, mark size= 0.2ex] (-3,0);
\node[align=left, below] at (-3,0){$0$};
\draw plot [mark=*, mark size= 0.2ex] (-1.5,0);
\node[align=left, below] at (-1.5,0){$l_1$};
\draw plot [mark=*, mark size= 0.2ex] (1.5,0);
\node[align=left, below] at (1.5,0){$l_2$};
\draw plot [mark=*, mark size= 0.2ex] (3,0);
\node[align=left, below] at (3,0){$i$};
\draw plot [mark=*, mark size= 0.2ex] (4,0);
\node[align=left, below] at (4,0){$l_4$};
\end{tikzpicture}
\caption{Graph of the return times to $I_{L+1}$}
\end{center}
\end{figure}
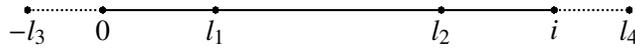

Let $\left\{{p_n}/{q_n}\right\}$ be the continued fraction approximation of $\alpha$. Assume that $q_L\leq i<q_{L+1}$.
By Lemma \ref{lemma.returntime} in the appendix, there are at most two $l$'s such that $0\leq l\leq i$ and $T^lx\in I_{L+1}$.
Let $l_1$, $l_2$ and $l_4$ be the first, second and third time for the trajectory entering $I_{L+1}$ forward. Let $l_3$ be the first time for the trajectory entering $I_{L+1}$ backward.
We focus on the case that there exist two such $l$'s, see Figure 1, and the other cases can be dealt with in a similar way.
Now the trajectory is divided into three segments:
$$\{T^jx:j\in[0,l_1]\},\quad \{T^jx:j\in[l_1, l_2]\},\quad \{T^jx:j\in[l_2, i]\}.$$
Accordingly, we consider the decomposition
\beq\label{composition}
A_i(x)=A_{i-l_2}(T^{l_2}x)A_{l_2-l_1}(T^{l_1}x)A_{l_1}(x).
\eeq
Let $s_n(x)=s(A_n(x))$ and $u_{n}(x)=s(A_{-n}(x))$. Using the polar decomposition, (\ref{composition}) becomes
\begin{small}
\beas
A_i(x)&=&R_{u_{i-l_2}(T^{i}x)}
\left(\begin{array}{cc}
    \|A_{i-l_2}\| & 0 \\
    0 & \|A_{i-l_2}\|^{-1} \\
  \end{array}\right)
R_{\frac\pi 2-s_{i-l_2}(T^{l_2}x)+u_{l_2-l_1}(T^{l_2}x)}\\
&&\cdot\left(\begin{array}{cc}
    \|A_{l_2-l_1}\| & 0 \\
    0 & \|A_{l_2-l_1}\|^{-1} \\
  \end{array}\right)
R_{\frac\pi 2-s_{l_2-l_1}(T^{l_1}x)+u_{l_1}(T^{l_1}x)}
\left(\begin{array}{cc}
    \|A_{l_1}\| & 0 \\
    0 & \|A_{l_1}\|^{-1} \\
  \end{array}\right)
R_{\frac\pi 2-s_{l_1}(x)}.
\eeas
\end{small}

\noindent Our basic idea is to show the norms of $A_{i-l_2}$, $A_{l_2-l_1}$ and $A_{l_1}$ are large enough and the corresponding angles are not too small and then apply Lemma \ref{lemma.basic}.
Depending on the positions of $l_1$ and $l_2$, we divide the proof into the following cases. See Figure 2.
\begin{description}
  \item [Case 1] $i-l_2,\ l_2-l_1, \ l_1 \ \geq 10^{-3} i$.
  \item [Case 2] $i-l_2, \ l_1 \ \geq 10^{-3} i, \quad l_2-l_1<10^{-3}i$.
  \item [Case 3] Other degenerate cases: $l_1<10^{-3} i$ or $i-l_2<10^{-3} i$; there is only one $l$'s such that $0\leq l\leq i$ and $T^lx\in I_{L+1}$; there is no $l$'s such that $0\leq l\leq i$ and $T^lx\in I_{L+1}$.
\end{description}
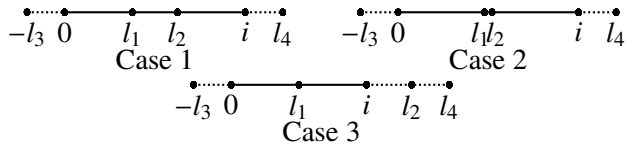
\begin{figure}[H]\label{lemma.LDT.case123}
\begin{center}
\begin{tikzpicture}[yscale=1.2]
\draw [densely dotted, thick] (-1.7,0) -- (-1.2,0);
\draw [ thick] (-1.2,0) -- (1.2,0);
\draw [densely dotted, thick] (1.2,0) -- (1.7,0);
\draw plot [mark=*, mark size= 0.2ex] (-1.7,0);
\node[align=left, below] at (-1.7,0){$-l_3$};
\draw plot [mark=*, mark size= 0.2ex] (-1.2,0);
\node[align=left, below] at (-1.2,0){$0$};
\draw plot [mark=*, mark size= 0.2ex] (-0.3,0);
\node[align=left, below] at (-0.3,0){$l_1$};
\draw plot [mark=*, mark size= 0.2ex] (0.3,0);
\node[align=left, below] at (0.3,0){$l_2$};
\draw plot [mark=*, mark size= 0.2ex] (1.2,0);
\node[align=left, below] at (1.2,0){$i$};
\draw plot [mark=*, mark size= 0.2ex] (1.7,0);
\node[align=left, below] at (1.7,0){$l_4$};
\node[align=left, below] at (0,-0.3){Case 1};
\end{tikzpicture}\ \ \ \
\begin{tikzpicture}[yscale=1.2]
\draw [densely dotted, thick] (-1.7,0) -- (-1.2,0);
\draw [ thick] (-1.2,0) -- (1.2,0);
\draw [densely dotted, thick] (1.2,0) -- (1.7,0);
\draw plot [mark=*, mark size= 0.2ex] (-1.7,0);
\node[align=left, below] at (-1.7,0){$-l_3$};
\draw plot [mark=*, mark size= 0.2ex] (-1.2,0);
\node[align=left, below] at (-1.2,0){$0$};
\draw plot [mark=*, mark size= 0.2ex] (-0.05,0);
\node[align=left, below] at (-0.1,0){$l_1$};
\draw plot [mark=*, mark size= 0.2ex] (0.05,0);
\node[align=left, below] at (0.1  ,0){$l_2$};
\draw plot [mark=*, mark size= 0.2ex] (1.2,0);
\node[align=left, below] at (1.2,0){$i$};
\draw plot [mark=*, mark size= 0.2ex] (1.7,0);
\node[align=left, below] at (1.7,0){$l_4$};
\node[align=left, below] at (0,-0.3){Case 2};
\end{tikzpicture}
\ \ \ \
\begin{tikzpicture}[yscale=1.2]
\draw [densely dotted, thick] (-1.7,0) -- (-1.2,0);
\draw [thick] (-1.2,0) -- (0.6,0);
\draw [densely dotted, thick] (0.6,0) -- (1.7,0);
\draw plot [mark=*, mark size= 0.2ex] (-1.7,0);
\node[align=left, below] at (-1.7,0){$-l_3$};
\draw plot [mark=*, mark size= 0.2ex] (-1.2,0);
\node[align=left, below] at (-1.2,0){$0$};
\draw plot [mark=*, mark size= 0.2ex] (-0.3,0);
\node[align=left, below] at (-0.3,0){$l_1$};
\draw plot [mark=*, mark size= 0.2ex] (0.6,0);
\node[align=left, below] at (0.6,0){$i$};
\draw plot [mark=*, mark size= 0.2ex] (1.2,0);
\node[align=left, below] at (1.2,0){$l_2$};
\draw plot [mark=*, mark size= 0.2ex] (1.7,0);
\node[align=left, below] at (1.7,0){$l_4$};
\node[align=left, below] at (0,-0.3){Case 3};
\end{tikzpicture}

\caption{The return times to $I_{L+1}$}
\label{f.graph-type123}
\end{center}
\end{figure}

\textbf{Case 1.}
Consider the decomposition
$$
A_{i}(x)=A_{i-l_2}(T^{l_2}x)A_{l_2-l_1}(T^{l_1}x)A_{l_1}(x).
$$
Since $l_2\geq 10^{-3}i\geq 10^{-3}q_L\geq r_L$, we have by Theorem \ref{theorem.induction} that
$$
%\|A_{l_1}(x)\|\geq \lambda^{\frac{99}{100}l_1}, \ \
%\|A_{i-l_2}(T^{l_2}x)\|\geq \lambda^{\frac{99}{100}(i-l_2)},\ \
\|A_{l_2-l_1}(T^{l_1}x)\|\geq \lambda^{\frac{99}{100}(l_2-l_1)}.
$$
We need the following lemma to obtain the estimate of $A_{l_1}(x)$. (Similarly for $A_{i-l_2}(T^{l_2}x)$.)
\begin{lemma}\label{keylemma}
Let $N$,$\lambda$ and $i$ be as in Theorem \ref{LDT}.
Let $L\in \mathbb Z$ satisfy $q_L\leq i <q_{L+1}$.
Assume that $l_1\geq 10^{-3} q_L$. Then we have
\bea\label{keylemma.norm}
&\|A_{l_1}(x)\|\geq \lambda^{\frac{99}{100}l_1},
\\\label{keylemma.angle}
&|u_{l_1}(T^{l_1}x)-u_{l_1+l_3}(T^{l_1}x)| <\lambda^{-\frac{99}{100}l_1}.
\eea
\end{lemma}
\noindent The proof of this lemma is so technical that we put it into Section 4.
\vskip 0.3cm
Next, we want to show the nondegenerate property of the angle between $A_{l_1}(x)$ and $A_{l_2-l_1}(T^{l_1}x)$.
Let
\beq\label{theorem.LDT.2}
f_1(T^{l_1}x)=s_{l_2-l_1}(T^{l_1}x)-u_{l_1+l_3}(T^{l_1}x).
\eeq
\begin{lemma}\label{lemma.keytypeI}
If $T^{l_1}x\in I_{L+1}\backslash B(C^{(L+2)}, \lambda^{-2\cdot 10^{-6}i})$, then
\beq\label{theorem.LDT.3}
|f_1(T^{l_1}x)|\geq C \lambda^{-2\cdot 10^{-6}i}.
\eeq
\end{lemma}
\begin{prof}
Since $l_2-l_1\geq 10^{-3}q_L$ and $l_1+l_3\geq 10^{-3}q_L$, then $l_2-l_1\geq r_L$ and $l_1+l_3\geq r_L$. Hence by Lemma \ref{lemma.basic} we have
$$
\|s_{l_2-l_1}(T^{l_1}x)-s_{r^+_{L}}(T^{l_1}x)\|_{C^2}\leq \lambda^{-\frac32 r_L},
$$
$$
\|u_{l_1+l_3}(T^{l_1}x)-u_{r^-_{L}}(T^{l_1}x)\|_{C^2}\leq \lambda^{-\frac32 r_L}.
$$
%Here we have $C^2$ estimate because ...
Then we have
\beq\label{lemma.keytypeI.1}
\|f_1(T^{l_1}x)-g_{L+1}(T^{l_1}x)\|_{C^2}\leq \lambda^{-\frac32 r_L}.
\eeq

Without loss of generality, we only focus on $I_{L+1,1}$.
According to \cite{WZ}, for each $j$, $g_{j}$ belongs to one of three cases, namely type $\I$, $\II$ and $\III$ (see appendix for definitions and properties of them).
For each case, we have the lower bound estimate for $|g_{j}(x)|$, $|\frac{dg_j}{dx}(x)|$ or $|\frac{d^2g_j}{dx^2}(x)|$.
If $g_{L+1}$ is of type $\I$ in $I_{L,1}$,
then $|\frac{dg_{L+1}}{dx}|\geq q_L^{-4}$. And hence $|\frac{df_1}{dx}|\geq q_L^{-5}$ by (\ref{lemma.keytypeI.1}). Then $f_1$ has only one zero, say $e_1$. If $T^{l_1}x\notin B(e_1, \lambda^{-2\cdot 10^{-6}i})$, then (\ref{theorem.LDT.3}) holds.

Now assume $g_{L+1}$ be of type $\III$ in $I_{L,1}$.
If $g_{L+2}$ is of type $\I$ in $I_{L+1,1}$, then by Lemma \ref{typeIII}, $g_{L+1}$ is of type $\I$ in $I_{L+1,1}$. This implies $|\frac{dg_{L+1}}{dx}|\geq q_{L+1}^{-4}\gg \lambda^{-\frac32 r_L}$ and hence $|\frac{df_1}{dx}|\geq q_{L+1}^{-5}$ by (\ref{lemma.keytypeI.1}). If $T^{l_1}x\notin B(e_1, \lambda^{-2\cdot 10^{-6}i})$, then (\ref{theorem.LDT.3}) holds.

If $g_{L+2}$ is also of type $\III$ in $I_{L+1,1}$.
Also by Lemma \ref{typeIII} in the appendix, $g_{L+1}$ is of type $\I$ in $I_{L+1,1}$, or of type $\III$ in $I_{L,1}$ with $\|g_{L+2}-g_{L+1}\|_{C^2}\leq \lambda^{-\frac32r_{L}}$.
On the other hand, since $l_2-l_1>10^{-3}q_L$, one can compare the derivatives of $g_{L+2}$ and $g_{L+1}$ to obtain  $\|g_{L+2}-g_{L+1}\|_{C^2}\not\leq \lambda^{-\frac32r_{L}}$. In fact, due to the form of type $\III$ function, $max_x|\frac{dg_{L+2}}{dx}|\approx \|A_{l_2-l_1}\|^{2} >\lambda^{10^{-3}q_L} \gg\lambda^{2\cdot10^{-6}q_L}>\max_x|\frac{dg_{L+2}}{dx}|.$
Hence $g_{L+1}$ is of type $\I$ in $I_{L+1,1}$, which is also reduced to the above case.

Now it is left to estimate the distance between $e_1$ and $c_{L+2,1}$.
Recall that $g_{L+2}$ has the following form:
$$
g_{L+2}(T^{l_1}x)= \arctan\left(\|A_{l_2-l_1}(T^{l_1}x)\|^2\tan f_2(T^{l_2}x)\right) -\frac \pi 2+f_1(T^{l_1}x).
$$
where $f_2(T^{l_2}x)=s_{l_4-l_2}(T^{l_2}x)-u_{l_2-l_1}(T^{l_2}x)$.
Using (\ref{lemma.typeIII.1}) in Lemma \ref{lemma.propertyoftpye3} in the appendix, we have
$$
|c_{L+2,1}-e_1|\leq \|A_{l_2-l_1}(T^{l_1}x)\|^{-\frac 34}<\lambda^{-\frac 34\cdot 10^{-3}i}\ll\lambda^{-2\cdot 10^{-6}i}.
$$
Hence if $T^{l_1}x\notin B(c_{L+2,1}, \lambda^{-2\cdot 10^{-6}i})$,
then (\ref{theorem.LDT.3}) holds.

\end{prof}

Combining (\ref{keylemma.angle}), (\ref{theorem.LDT.2}) and (\ref{theorem.LDT.3}), it holds that
$$
|s_{l_2-l_1}(T^{l_1}x)-u_{l_1}(T^{l_1}x)|\geq C\lambda^{-2\cdot 10^{-6}i}-\lambda^{-\frac{99}{100}\cdot10^{-3}i}\geq C\lambda^{-2\cdot 10^{-6}i}.
$$
Similarly, if $T^{l_2}x\notin B(c_{L+2,2}, \lambda^{-2\cdot 10^{-6}i})$,
then $|f_2(T^{l_2}x)|\geq \lambda^{-2\cdot 10^{-6}i}$ and
$$|s_{i-l_2}(T^{l_2}x)-u_{l_2-l_1}(T^{l_2}x)|\geq C\lambda^{-2\cdot 10^{-6}i}.$$
Therefore if
$$x\notin  \Omega_1:=\bigcup_{l_1,l_2\in[0,i]}\left(B(C^{(L+2)},\lambda^{-2\cdot 10^{-6}i})-l_1\alpha\right) \cup\left(B(C^{(L+2)},\lambda^{-2\cdot 10^{-6}i})-l_2\alpha\right),$$
then one can apply the nonresonant case of Lemma \ref{lemma.basic} to obtain
$$
\left\|A_{i}(x)\right\|\geq \lambda^{\frac9{10}i}.
$$
The measure of the exception set is $$\text{meas} \ \Omega_1 <i^2\lambda^{-2\cdot 10^{-6}i} <\lambda^{-\delta i}.$$
This implies the result.

\

\textbf{Case 2.}
Since $i-l_2\geq \frac12i \geq 10^{-3}i\geq l_2-l_1$, then
$$\| A_{i-l_2}(T^{l_2}x)\|\gg \| A_{l_2-l_1}(T^{l_1}x)\|.$$
It follows from Lemma \ref{keylemma} and Lemma \ref{lemma.basic} that
\beqs
\| A_{i-l_1}(T^{l_1}x) \| \geq \lambda^{\frac{98}{100}(i-l_1)}\geq \lambda^{\frac{49}{100}i}.
\eeqs
\beq\label{theorem.LDT.4}
|u_{i-l_2}(T^{i}x)-u_{i-l_1}(T^{i}x)| < \|A_{i-l_2}(T^{i}x)\|^{-1}<\lambda^{-\frac12 i}.
\eeq
Later we need
\begin{lemma}\label{lemma.LDT.5}
\beq\label{theorem.LDT.5}
|s_{l_4-i}(T^{i}x)-u_{i-l_2}(T^{i}x)| \geq C\left\|T^{i}x-C^{(i_n)}\right\|\geq Cq_{i_{n+1}}^{-2}.
\eeq
\end{lemma}
\noindent We put the proof at the end of section 4.
\vskip 0.5cm
Combining (\ref{theorem.LDT.4}) and Lemma \ref{lemma.LDT.5},
$$
|s_{l_4-i}(T^{i}x)-u_{i-l_1}(T^{i}x)|\geq Cq_{i_{n+1}}^{-2}.
$$
Since $l_4-i\geq 10^{-3}i \geq 10^{-3}q_L$, by Lemma \ref{keylemma} we have that $\| A_{l_4-i}(T^ix)\| \geq \lambda^{\frac{99}{100}(l_4-i)}$.
Then it follows by Lemma \ref{lemma.basic} that
$$
|s_{l_4-l_1}(T^{l_1}x)-s_{i-l_1}(T^{l_1}x)| <\|A_{i-l_1}(T^{l_1}x)\|^{-\frac32} \leq \lambda^{-\frac{49}{100}i}.
$$
We also obtain by Lemma \ref{keylemma} that
$$
|u_{l_1}(T^{l_1}x)-u_{l_1+l_3}(T^{l_1}x)| < \lambda^{-\frac{99}{100} \cdot10^{-3}i}.
$$
Then the above two inequalities imply that
$$
|f_3(T^{l_1}x)-g_{L+2}(T^{l_1}x)| < \lambda^{-\frac{99}{100} \cdot10^{-3}i}.
$$
where $f_3(T^{l_1}x)=s_{i-l_1}(T^{l_1}x)-u_{l_1}(T^{l_1}x)$. In this case, $g_{L+2}$ is of type $\III$.
Let $$ x\notin \Omega_3:=\bigcup_{ l_1,l_2\in[0,i]} \left(B(C^{(L+2)},\lambda^{-2\cdot 10^{-6}i}) -l_1\alpha\right). $$
Then
$$
|f_3(T^{l_1}x)|\geq |g_{L+2}(T^{l_1}x)|-\lambda^{-\frac{99}{100} \cdot10^{-3}i} \geq C\lambda^{-6\cdot10^{-6}i}-\lambda^{-\frac{99}{100} \cdot10^{-3}i} \geq C\lambda^{-6\cdot10^{-6}i}.
$$
Then we obtain
$$
\|A_{i}(x)\|\geq \|A_{i-l_1}(T^{l_1}x)\|\cdot \|A_{l_1}(x)\| \cdot \lambda^{-6\cdot10^{-6}i}
\geq \lambda^{\frac{9}{10}i}.
$$
The measure of the exception set is
$$
\text{meas} \ \Omega_3<i^2\lambda^{-2\cdot 10^{-6}i}<\lambda^{-\delta i}.
$$

\

\textbf{Case 3.}
Actually this case is easier than Case 1.
We first focus on the case that there is only one point of the trajectory entering $I_{L+1}$ in $[0,i]$, and it holds that $$ l_1\geq 10^{-3}i,\ \ i-l_1\geq 10^{-3}i. $$
See Figure 2. We follow the same idea of Case 1.
By Lemma \ref{keylemma}, it holds that
\beas
|u_{l_1}(T^{l_1}x)-u_{l_1+l_3}(T^{l_1}x)| <\lambda^{-\frac{99}{100}\cdot10^{-3}i},
\\
|s_{i-l_1}(T^{l_1}x)-s_{l_2-l_1}(T^{l_1}x)| <\lambda^{-\frac{99}{100}\cdot10^{-3}i}.
\eeas
Let $$f_4(T^{l_1}x)=s_{l_2-l_1}(T^{l_1}x)-u_{l_1+l_3}(T^{l_1}x).$$
By Lemma \ref{lemma.keytypeI},
if $T^{l_1}x\notin B(c_{L+2,1},\lambda^{-2\cdot 10^{-6}i})$, then
$|f_4(T^{l_1}x)|\geq \lambda^{-2\cdot 10^{-6}i}$.
This implies
$$
|s_{i-l_1}(T^{l_1}x)-u_{l_1}(T^{l_1}x)|\geq C\lambda^{-2\cdot 10^{-6}i}-\lambda^{-\frac{99}{100}\cdot10^{-3}i}\geq C\lambda^{-2\cdot 10^{-6}i}.
$$
By Lemma \ref{lemma.basic}, we conclude that if $$x\notin \Omega_2:= \bigcup_{l_1\in[0,i]}\left(B(C^{(L+2)},\lambda^{-2\cdot 10^{-6}i})-l_1\alpha\right),$$ then
$$
\|A_{i}(x)\|\geq \lambda^{\frac9{10}i}.
$$
The measure of the exception set is $$\text{meas} \ \Omega_2 <i\lambda^{-2\cdot 10^{-6}i} <\lambda^{-\delta i}.$$
Therefore the result follows.

For the case with $l_1<10^{-3} i$ or $i-l_2<10^{-3} i$, one can follow the same approach above to obtain the result.
For the case that there is no $l$'s such that $0\leq l\leq i$ and $T^lx\in I_{L+1}$, we can apply Lemma \ref{lemma.basic} to obtain the norm of $A_i(x)$ with no exception set of $x$ removed.
\end{prof}

%\begin{remark}\label{Holder exponent}
%It is easy to see that the value of the coefficient $C=10^{-6}$ is related to the coefficient $\kappa$ in Lemma \ref{lemma.basic}, see Remark \ref{kappa}.
%\end{remark}

\section{Proof of Lemma \ref{keylemma} and \ref{lemma.LDT.5}}

Let us first outline the proof of Lemma \ref{keylemma}.
Denote by $j_k$ $(k\leq L+1)$ the first time for the trajectory $\{T^jx:j\in \mathbb Z\}$ entering $I_k$.
Note that if $j_{k+1}-j_{k}\geq r_{k}^+(T^{j_k}x)$, then we can apply Theorem \ref{theorem.induction} to obtain the estimate of $\|A_{j_{k+1}-j_{k}}\|$, $s_{j_{k+1}-j_{k}}$ and $u_{j_{k+1}-j_{k}}$.

The difficulty arises when $j_{k+1}-j_{k}<r_{k}^+(T^{j_k}x)$. The condition of Theorem \ref{theorem.induction} fails due to the occurrence of the resonance in the $q_{k}$'s scale.
Since the frequency is Liouvillean, we have to analyze all these resonances, which is very complicated.
Inspired by the technique of "CD-bridge" in \cite{AFK}, we simplify the procedure by defining $i_l$ ($l\leq n+1$) as
$$q_{i_l}=\left\{
          \begin{array}{ll}
            q_{L+1}, & \hbox{$l=n+1$;} \\
            \max\left\{q_k: q_k^3\leq q_{i_{l+1}}\right\}, & \hbox{$l\leq n$.}
          \end{array}
        \right.
$$
With this definition,
\begin{itemize}
  \item It holds that $j_{i_{n-2}}\leq i^{\frac13}$. Hence $j_{i_{n-2}}$ is small enough compared with $i$ and then without loss of generality in the following we assume that $x\in I_{i_{n-2}}$;
  \item We can well estimate $A_{j_{i_{l+1}}-j_{i_l}}$ as the following
      $$\|A_{j_{i_{l+1}}-j_{i_l}}\|\gtrsim\lambda^{\frac{99}{100}(j_{i_{l+1}}-j_{i_l})},\quad s_{j_{i_{l+1}}-j_{i_l}}\thickapprox s_{r^+_{i_l}}, \quad u_{j_{i_{l+1}}-j_{i_l}}\thickapprox u_{r^-_{i_l}}.$$
\end{itemize}
Therefore we only need to consider the decomposition
%\begin{tiny}
\bea\nn
A_{l_1+l_3}(T^{-l_3}x)&=&A_{l_1}(x)A_{l_3}(T^{-l_3}x)\\
\label{keylemma.decomposition}&=&A_{l_1-j_{i_n}}(T^{j_{i_n}}x)A_{j_{i_n}-j_{i_{n-1}}}(T^{j_{i_{n-1}}}x) A_{j_{i_{n-1}}}(x) \\
&&\ \cdot \ A_{j_{-i_{n-1}}}(T^{j_{-i_{n-1}}}x)A_{j_{-i_n}-j_{-i_{n-1}}}(T^{j_{-i_{n}}}x)A_{l_3-j_{-i_n}}(T^{-l_3}x).\nn
\eea
%\end{tiny}
where $j_{-k}$ $(k\leq L+1)$ is the first time for the trajectory $\{T^jx:j\in \mathbb Z\}$  entering $I_k$ backward.
Note if we have $|s_{l_1}(x)-u_{l_3}(x)|\geq \max\left\{\|A_{l_1}(x)\|^{-\eta}, \|A_{l_3}(T^{-l_3}x)\|^{-\eta}\right\}$,
then the result holds by Lemma \ref{lemma.basic}.
For this purpose, we will repeatedly apply Lemma \ref{lemma.basic} to show in the nondegenerate case
$$
s_{l_1}(x)\approx s_{j_{i_n}}(x)\approx s_{j_{i_{n-1}}}(x),\quad
u_{l_3}(x)\approx u_{j_{-i_n}}(x)\approx u_{j_{-i_{n-1}}}(x).
$$
%In this procedure, one can also easily obtain the estimate (\ref{keylemma.norm}).
On the other hand, by Theorem \ref{theorem.induction},
$$
|s_{j_{i_{n-1}}}(x)-u_{j_{-i_{n-1}}}(x)|\geq \max\left\{\|A_{l_1}(x)\|^{-\eta}, \|A_{l_3}(T^{-l_3}x)\|^{-\eta}\right\}.
$$
This finishes the proof.
We remark at last that if the frequency is Diophantine, then the proof of Lemma \ref{keylemma} can be greatly simplified.

\

Now we turn to the formal proof.
Let $\lambda_1:=\lambda_3(\alpha,v ,10^{-3})$ and $\lambda>\lambda_1$.
Let  $L$ satisfy  $q_L\leq i<q_{L+1}$.
It holds that $q_{i_{l+1}}\geq q_{i_{l}}^3$ for $ l\leq n$. The following lemma leads to the upper bound of the return time.
\begin{lemma}\label{lemma.qin}
Assume that $x$ is an arbitrary point on $\R/\Z$ and $I=B\left(c,q_{i_l}^{-2}\right)$. Then we have
$$ \min\left\{j\in \mathbb N:T^jx\in I\right\}<q_{i_{l+1}}.$$
\end{lemma}

\begin{prof}
According to (\ref{arithmetical.1}) in the appendix, we have
$$
\|k\alpha-\frac{kp_{i_{l+1}}}{q_{i_{l+1}}}\| <\frac{k}{q_{i_{l+1}}q_{i_{l+1}+1}}<\frac{1}{q_{i_{l+1}+1}},\quad 0<k<q_{i_{l+1}}.$$
Since ${q_{i_{l+1}}^{-1}}\leq {q_{i_{l}}^{-3}} \ll {q_{i_{l}}^{-2}}$, then there exists $k_0$ such that $0<k_0<q_{i_{l+1}}$ and
$$
\left(\frac{k_0p_{i_{l+1}}}{q_{i_{l+1}}}-\frac{1}{q_{i_{l+1}+1}}, \frac{k_0p_{i_{l+1}}}{q_{i_{l+1}}}+\frac{1}{q_{i_{l+1}+1}}\right)\subset I.
$$
The lemma is now evident from what we have proved.
\end{prof}

By Lemma \ref{lemma.qin}, it is straightforward to see that
$$
\min\left\{j\in \mathbb N: T^{j}x\in I_{{i_{n-2}}}\right\}\leq q_{i_{n-1}}\leq q_{i_n}^\frac13 \leq q_L^\frac13 \leq i^{\frac13}.
$$
And hence $$\log \|A_{j_{i_{n-2}}}\|\leq i^{\frac13} \log \lambda \ll i. $$
In other words, $\|A_{j_{i_{n-2}}}\|$ is small enough.
Then we can assume $x\in I_{i_{n-2}}$.
Let $j_k$ be the first time for the trajectory of the base dynamic entering $I_k$.
Then it holds that $j_{i_{n-2}}=0$ and $j_{i_{n+1}}=l_1$.
Hence we only need to consider the decomposition (\ref{keylemma.decomposition}).

By the definition of $i_l$,
$$ q_{i_{l+1}}\leq q_{i_l+1}^3 \leq \exp(3\cdot 10^{-12}\log\lambda \cdot q_{i_l}).$$
The second inequality requires
\beqs
\log\lambda \geq 10^{12} q_{i_{l}}^{-1}\log q_{i_{l}+1},
\eeqs
which can be derived from (\ref{thm.lambda0}) and Remark \ref{2.5}.
Then
$$
i< q_{i_{n+1}}\leq \exp(3\cdot 10^{-12}\log\lambda \cdot \exp(3\cdot 10^{-12}\log\lambda\cdot\exp(3\cdot 10^{-12}\log\lambda \cdot q_{i_{n-2}}))).
$$
Let $$i_0=\exp(3\cdot 10^{-12}\log\lambda\cdot \exp(3\cdot 10^{-12}\log\lambda\cdot\exp(3\cdot 10^{-12}\log\lambda \cdot q_{i_{n-2}}))).$$
For $i>i_0$, we have $q_{i_{n-2}}\geq q_N$.
Then Theorem \ref{theorem.induction} is applicable for $I_{i_{n-2}}$.

\begin{lemma} \label{claim.1}
If $j_{i_{n-1}}\geq 10^{-6}q_{i_{n-2}}$, then it holds that
\bea \label{claim.1.4}
&\|A_{j_{i_{n-1}}}(x)\| \geq\lambda^{(1-10^{-4})j_{i_{n-1}}},
\\ \label{claim.1.5}
&|s_{j_{i_{n-1}}}(x) -s_{r_{i_{n-2}}^+}(x)|<\lambda^{-\frac32 r_{i_{n-2}}}.
%\\ \label{claim.1.6} &|u_{j_{i_{n-1}}}(T^{j_{i_{n-1}}}x) - u_{r_{i_{n-2}}^-}(T^{j_{i_{n-1}}}x)|<\lambda^{-\frac32 r_{i_{n-2}}}.
\eea
Moreover, similar estimates hold for $A_{j_{i_{n}}-j_{i_{n-1}}}(T^{j_{i_{n-1}}}x)$ if $j_{i_{n}}-j_{i_{n-1}}\geq 10^{-6} q_{i_{n-1}}$ and for $A_{l_1-j_{i_n}}(T^{j_{i_{n}}}x)$ if $l_1-j_{i_n}\geq 10^{-6}q_{i_n}$.
\end{lemma}
\begin{prof}
By the definition of $j_{i_{n-1}}$, we have $T^jx\notin I_{i_{n-1}}$ for $0<j<j_{i_{n-1}}$, which implies
\beq\label{claim.1.0}
\left\|T^{j}x-C^{(i_{n-1})}\right\|\geq q_{i_{n-1}}^{-2}\geq q_{i_{n-2}+1}^{-6}\geq \lambda^{-6\cdot 10^{-12}q_{i_{n-2}}}.
\eeq
Since $j_{i_{n-1}}\geq 10^{-6}q_{i_{n-2}}$, then $j_{i_{n-1}}\geq r_{i_{n-2}}^+(x)$.
Let $$0=t_0<t_1<t_2<\cdots<t_p\leq j_{i_{n-1}}$$ be the return times to $I_{i_{n-2}}$ such that $t_{l+1}-t_{l}\geq10^{-6}{q_{i_{n-2}}}$ ($0\leq l <p$) and $0\leq j_{i_{n-1}}-t_p<10^{-6}{q_{i_{n-2}}}$.
\begin{itemize}
  \item
  If $t_p= j_{i_{n-1}}$,
  then by Theorem \ref{theorem.induction}, we obtain
  \bea\label{claim.1.1}
  &\|A_{t_{l}-t_{l-1}}(T^{t_{l-1}}x)\|\geq \lambda^{(1-10^{-5})(t_{l}-t_{l-1})},\\
  \label{claim.1.2}
  &|s_{t_{l+1}-t_{l}}(T^{t_{l}}x)-u_{t_{l}-t_{l-1}}(T^{t_{l}}x) -g_{i_{{n-1}}+1}(T^{t_{l}}x)| \leq \lambda^{-\frac32r_{i_{n-2}}}.
  \eea
  And hence
  \bea\label{claim.1.3}
 |s_{t_{l+1}-t_{l}}(T^{t_{l}}x)-u_{t_{l}-t_{l-1}}(T^{t_{l}}x))|&\geq & |g_{i_{{n-1}}+1}(T^{t_{l}}x)|-\lambda^{-\frac32r_{i_{n-2}}}\\ &\geq&
  C q_{i_{n-1}}^{-6}-\lambda^{-\frac32 \cdot 10^{-6}{q_{i_{n-2}}} }\nn\\ &\geq& C q_{i_{n-1}}^{-6}.\nn
  \eea
  In the second inequality we use (\ref{nondegenerate}) and (\ref{claim.1.0}).
  Then one can easily verify the conditions of the nonresonant case in Lemma \ref{lemma.basic} and the result follows.

  \item If $t_p<j_{i_{n-1}}$. It holds that
      $$j_{i_{n-1}}-10^{-6}{q_{i_{n-2}}}<t_p,\ \  q_{i_{n-2}}-10^{-6}{q_{i_{n-2}}}\leq t_p.$$
      It is easy to see that (\ref{claim.1.1}), (\ref{claim.1.2}) and (\ref{claim.1.3}) still hold. Then we can apply Lemma \ref{lemma.basic} to obtain
      \beqs
      \|A_{t_p}(x)\| \geq\lambda^{(1-2\cdot10^{-5})t_p},\quad
     |s_{t_p}(x) -s_{r_{i_{n-2}}^+}(x)|<\lambda^{-\frac32 r_{i_{n-2}}}.
      \eeqs
      Since $j_{i_{n-1}}-t_p<10^{-6}q_{i_{n-2}}< 10^{-5}t_p$, we obtain (\ref{claim.1.4}).
      Moreover, by Lemma \ref{lemma.basic},
      $$
      |s_{j_{i_{n-1}}}(x) -s_{t_p}(x)|<\lambda^{-\frac32 r_{i_{n-2}}}.
      $$
      And hence (\ref{claim.1.5}) follows.
      %To obtain (\ref{claim.1.6}), one can consider $A_{-j_{i_{n-1}}}(T^{j_{i_{n-1}}}x)$ and the trajectory $\{T^{j}x:0\leq j \leq j_{i_{n-1}}\}$ backwards. Repeat the above process, (\ref{claim.1.6}) follows as the similar way to prove (\ref{claim.1.5}).
\end{itemize}
\end{prof}

From the proof of Lemma \ref{claim.1}, one can immediately obtain the following.
\begin{corollary}\label{claim.3}
  Assume that $x\in I_{i_{n-2}}$, $q_{i_{n}+1}\leq q_{i_{n}}^3$ and $l_1-j_{i_{n-1}}\geq 10^{-6}q_{i_{n-1}}$. Then it holds that
  \beas
  &\|A_{l_1-j_{i_{n-1}}}(T^{j_{i_{n-1}}}x)\| \geq \lambda^{(1-10^{-4})(l_1-j_{i_{n-1}})},\\
  &|s_{l_1-j_{i_{n-1}}}(T^{j_{i_{n-1}}}x) -s_{r^+_{i_{n-1}}}(T^{j_{i_{n-1}}}x)| < \lambda^{-\frac32 r_{i_{n-1}}}.
  \eeas
  Similarly, if $q_{i_{n-1}+1}\leq q_{i_{n-1}}^3$ and  $j_{i_{n}}\geq 10^{-6}q_{i_{n-2}}$, then we have
  \beas
  &\|A_{j_{i_{n}}}(x)\| \geq \lambda^{(1-10^{-4})j_{i_{n}}},\\
 &| s_{j_{i_{n}}}(x) -s_{r^+_{i_{n-2}}}(x)| < \lambda^{-\frac32 r_{i_{n-1}}}.
  \eeas
\end{corollary}
\begin{prof}
Since $q_{i_{n}+1}\leq q_{i_{n}}^3$, by the definition of $i_l$, it holds that $q_{i_{n+1}}\leq q_{i_{n}+1}^3\leq q_{i_{n}}^9$.
Then $T^jx\notin I_{i_{n+1}}$ for $j_{i_{n-1}}\leq j< l_1$. This implies that
$$
\|T^{j}x-C^{(i_{n+1})}\|\geq q_{i_{n+1}}^{-2}\geq q_{i_{n}}^{-18} \geq q_{i_{n-1}+1}^{-54}\geq \lambda^{-54\cdot 10^{-12}q_{i_{n-1}}}.
$$
Using the same argument in the proof of Lemma \ref{claim.1}, this corollary follows.
\end{prof}

Corollary \ref{claim.3} shows that
\begin{itemize}
  \item if $q_{i_{n}+1}\leq q_{i_{n}}^3$, we only need to consider the decomposition $$A_{l_1}(x)=A_{l_1-j_{i_{n-1}}}(T^{j_{i_{n-1}}}x)A_{j_{i_{n-1}}}(x); $$
  \item if $q_{i_{n-1}+1}\leq q_{i_{n-1}}^3$, we only need to consider the decomposition $$A_{l_1}(x)=A_{l_1-j_{i_{n}}}(T^{j_{i_{n}}}x)A_{j_{i_{n}}}(x). $$
\end{itemize}
This consideration may simplify our proof.
From this we can see if the frequency is a Diophantine number, then the proof of Lemma \ref{keylemma} can be greatly simplified.

\

We are ready to prove the norm estimate (\ref{keylemma.norm}) and the angle estimate (\ref{keylemma.angle}).
We will apply Lemma \ref{lemma.basic} repeatedly to obtain these estimates.
At first, we would like to get the result in a special case.
%Denote by $j_{-k}$ the return time for the trajectory entering $I_k$ backwards.

\begin{lemma}\label{lemma.specialcase}
If $x\in I_{i_{n-1}}$, then the estimates (\ref{keylemma.norm}) and (\ref{keylemma.angle}) hold.
\end{lemma}
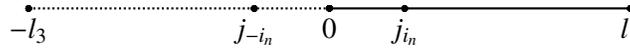
\begin{figure}[H]
\begin{center}
\begin{tikzpicture}[yscale=1]
\draw [densely dotted, thick] (-4,0) -- (0,0);
\draw [thick] (0,0) -- (4,0);
\draw plot [mark=*, mark size= 0.2ex] (-4,0);
\node[align=left, below] at (-4,0){$-l_3$};
\draw plot [mark=*, mark size= 0.2ex] (0,0);
\node[align=left, below] at (0,0){$0$};
\draw plot [mark=*, mark size= 0.2ex] (4,0);
\node[align=left, below] at (4,0){$l_1$};
\draw plot [mark=*, mark size= 0.2ex] (-1,0);
\node[align=left, below] at (-1,0){$j_{-i_{n}}$};
\draw plot [mark=*, mark size= 0.2ex] (1,0);
\node[align=left, below] at (1,0){$j_{i_{n}}$};
\end{tikzpicture}
\caption{Graph of the return times}
\end{center}
\end{figure}
\begin{prof}
We divide the proof into three cases: $x\in I_{i_{n-1}}\backslash I_{i_{n}}$, $x\in I_{i_n}\backslash I_{i_{n+1}}$ and $x\in I_{i_{n+1}}$.

If $x\in I_{i_{n+1}}$, then by definition we have $l_1=j_{i_{n+1}}=0$. There is nothing to say.

If $x\in I_{i_n}\backslash I_{i_{n+1}}$, then it follows that
$$
\left\|x-C^{(i_{n+1})}\right\|\geq q_{i_{n+1}}^{-2}\geq q_{i_{n}+1}^{-6}\geq \lambda^{-6\cdot 10^{-12} q_{i_n}}.
$$
If $l_1\leq10^{-3}q_{i_n}$ or $l_3\leq 10^{-3}q_{i_n}$, then $\|A_{l_1}(x)\|$ or $\|A_{i-l_3}(T^{l_3}x)\|$ is small enough and we can further assume $x\in I_{i_{n+1}}$, which is reduced to the previous case.
If $l_1 > 10^{-3} q_{i_n}$ and $l_3>10^{-3}q_{i_n}$, then by Lemma \ref{claim.1}, it holds that
\beas
&\| A_{l_1}(x) \| \geq \lambda^{(1-10^{-4})l_1},\ \ \| A_{-l_3}(x) \| \geq \lambda^{(1-10^{-4})l_3} ,
\\
&|s_{l_1}(x)-s_{r^+_{i_n}}(x)|<\lambda^{-\frac32r_{i_n}},
\\
&|u_{l_3}(x)-u_{r^-_{i_n}}(x)|<\lambda^{-\frac32r_{i_n}}.
\eeas
Hence
\beq\label{lemma.specialcase.0}
|s_{l_1}(x)-u_{l_3}(x)|\geq |g_{i_n+1}(x)|-\lambda^{-\frac32r_{i_n}}
\geq Cq_{i_{n+1}}^{-6}.
\eeq
By Lemma \ref{lemma.basic}, the estimate (\ref{keylemma.angle}) follows.

If $x\in I_{i_{n-1}}\backslash I_{i_{n}}$, then we have
$$
\left\|x-C^{(i_n)}\right\|\geq q_{i_n}^{-2}\geq q_{i_{n-1}+1}^{-6} \geq \lambda^{-6\cdot 10^{-12} q_{i_{n-1}}}.
$$
If $j_{i_n}\leq 10^{-3}q_{i_n}$ or $j_{-i_n}\leq 10^{-3}q_{i_n}$, then we can further assume that $x\in I_{i_n}\backslash I_{i_{n+1}}$, which we have already dealt with in the above.
Otherwise we have $j_{i_n}> 10^{-3}q_{i_n}$ and $j_{-i_n}> 10^{-3}q_{i_n}$, see Figure 3.
Again by Lemma \ref{claim.1}, it holds that
\bea\label{lemma.specialcase.1}
&\| A_{j_{i_n}}(x) \| \geq \lambda^{(1-10^{-4})j_{i_n}},\ \ \| A_{-j_{-i_n}}(x) \| \geq \lambda^{(1-10^{-4})j_{-i_n}} ,
\\\label{lemma.specialcase.2}
&|s_{j_{i_n}}(x)-s_{r^+_{i_{n-1}}}(x)|<\lambda^{-\frac32r_{i_{n-1}}},
\\\label{lemma.specialcase.3}
&|u_{j_{-i_n}}(x)-u_{r^-_{i_{n-1}}}(x)|<\lambda^{-\frac32r_{i_{n-1}}}.
\eea
And hence it holds that
$$
|s_{j_{i_n}}(x)-u_{j_{-i_n}}(x)| \geq |g_{i_{n-1}+1}(x)|-\lambda^{-\frac32r_{i_{n-1}}} \geq Cq_{i_{n}}^{-6}.
$$
Therefore one can apply Lemma \ref{lemma.basic} to obtain
\beq\label{lemma.specialcase.7}
|u_{j_{i_{n}}}(T^{j_{i_n}}x)-u_{r^-_{i_{n}}}(T^{j_{i_n}}x)| <\|A_{j_{i_n}}(x)\|^{-\frac32}<\lambda^{-10^{-3}q_{i_n}}.
\eeq
If $l_1-j_{i_n}<10^{-6}q_{i_n}$, then $$\|A_{l_1-j_{i_n}}(T^{j_{i_{n}}}x)\| \leq \lambda^{10^{-6}q_{i_n}} \ll \lambda^{(1-10^{-4})j_{i_n}}\leq \|A_{j_{i_n}}(x)\|.$$
Hence by Lemma \ref{lemma.basic}, it holds that
\bea
&\| A_{l_1}(x) \| \geq \|A_{j_{i_n}}\|\cdot\|A_{l_1-j_{i_n}}(T^{j_{i_{n}}}x)\|^{-1}\geq\lambda^{(1-10^{-4})j_{i_n}-10^{-6}q_{i_n}} \geq \lambda^{\frac{99}{100}l_1},
\nn\\
\label{lemma.specialcase.4}
&|s_{j_{i_n}}(x)-s_{l_1}(x)|<\|A_{j_{i_n}}(x)\|^{-\frac32}<\lambda^{-10^{-6}q_{i_n}},
\eea
where the first estimate is (\ref{keylemma.norm}).
Otherwise, we have $l_1-j_{i_n}\geq 10^{-6}q_{i_n}$, which by Lemma \ref{claim.1} implies that
\beq\label{lemma.specialcase.8}
|s_{l_1-j_{i_n}}(T^{j_{i_n}}x) - s_{r^+_{i_n}}(T^{j_{i_n}}x)| <\lambda^{-\frac32r_{i_n}}.
\eeq
Since $T^{j_{i_n}}x\in I_{i_n}\backslash I_{i_{n+1}}$, then $\|T^{j_{i_n}}x-C^{(i_{n+1})}\|\geq q_{i_{n+1}}^{-2}$. Combining (\ref{lemma.specialcase.7}) and (\ref{lemma.specialcase.8}), we obtain
\beas
&&|s_{l_1-j_{i_n}}(T^{j_{i_n}}x)-u_{j_{i_n}}(T^{j_{i_n}}x))|\\
&\geq& |g_{i_n+1}(T^{j_{i_n}}x)|-\lambda^{-\frac32r_{i_n}}-\lambda^{- 10^{-6}q_{i_n}}\\
&\geq& Cq_{i_{n+1}}^{-6}.
\eeas
By (\ref{lemma.specialcase.1}) and Lemma \ref{claim.1}, we have
\beqs
\|A_{j_{i_n}}(x)\|\geq \lambda^{(1-10^{-4})j_{i_n}},\ \ \|A_{l_1-j_{i_n}}(x)\|\geq \lambda^{(1-10^{-4})(l_1-j_{i_n})}.
\eeqs
Again one can apply Lemma \ref{lemma.basic} to obtain (\ref{keylemma.norm}) and (\ref{lemma.specialcase.4}). Similarly, we obtain
\bea
&\| A_{-l_3}(x) \| \geq \lambda^{\frac{99}{100}l_3},
\nn\\\label{lemma.specialcase.5}
&|u_{j_{-i_{n}}}(x)-u_{l_3}(x)|<\lambda^{-10^{-6}q_{i_n}}.
\eea
Combining (\ref{lemma.specialcase.2}), (\ref{lemma.specialcase.3}), (\ref{lemma.specialcase.4}) and (\ref{lemma.specialcase.5}), we have that
\bea\label{lemma.specialcase.6}
|s_{l_1}(x)-u_{l_3}(x)| &\geq&  |s_{j_{i_n}}(x)-u_{j_{-i_n}}(x)|-\lambda^{-10^{-6}q_{i_n}}\\\nn
&\geq& |g_{i_{n-1}+1}(x)|-\lambda^{-\frac32 r_{i_{n-2}}}-2\lambda^{- 10^{-6}q_{i_n}}\\
&\geq& C q_{i_{n}}^{-6}.\nn
\eea
Hence by Lemma \ref{lemma.basic} the estimate (\ref{keylemma.angle}) follows.
\end{prof}

\

We are now in a position to prove (\ref{keylemma.norm}) and (\ref{keylemma.angle}) in the case with $x\in I_{i_{n-2}}\backslash I_{i_{n-1}}$.
First we consider the case with $q_{i_{n}+1}> q_{i_{n}}^3$ and $q_{i_{n-1}+1}> q_{i_{n-1}}^3$. We divide the proof into the following cases.

\begin{enumerate}
  \item If $j_{i_{n-1}}< 10^{-3}q_{i_n}$ or $j_{-i_{n-1}}< 10^{-3}q_{i_n}$, then it holds that
      $$ \|A_{j_{i_{n-1}}}(x)\| \leq \lambda^{j_{i_{n-1}}} \leq \lambda^{10^{-3}q_{i_n}}\ \ \text{or} \ \ \|A_{-j_{-i_{n-1}}}(x)\| \leq \lambda^{j_{-i_{n-1}}}\leq \lambda^{10^{-3}q_{i_n}}. $$
      Hence we can reduce it to the case with $x\in I_{i_{n-1}}$, which has been dealt with in Lemma \ref{lemma.specialcase}.

  \item If $j_{i_n}-j_{i_{n-1}}< q_{i_{n-1}}$, then we claim $j_{i_n}\leq q_{i_{n-1}}$, which implies that
       $$\|A_{j_{i_n}}\|\leq \lambda^{q_{i_{n-1}}}\ll \lambda^{10^{-3}q_{i_n}}.$$
       Hence it can be reduced to the case with $x\in I_{i_{n}}$, from which we can obtain the estimates (\ref{keylemma.norm}) and (\ref{keylemma.angle}) by Lemma \ref{lemma.specialcase}.
      Suppose the claim was not true, then it holds that $j_{i_{n}}-q_{i_{n-1}}>0$.
      By the definition of $j_{i_n}$, $T^{j_{i_n}}x\in I_{i_n}$.
      Then we have
      \beas
      \left\|T^{j_{i_n}-q_{i_{n-1}}}x-C^{(i_{n-1})}\right\| &\leq& \left\|T^{j_{i_n}}x-q_{i_{n-1}}\alpha-C^{(i_{n})}\right\| +\left\|C^{(i_n)}-C^{(i_{n-1})}\right\|\\
      &\leq& \left\|T^{j_{i_n}}x-C^{(i_n)}\right\| +\left\|q_{i_{n-1}}\alpha\right\|+\lambda^{-\frac34r_{i_{n-1}}}\\
      &\leq& q_{i_n}^{-2}+q_{i_{n-1}+1}^{-1}+\lambda^{-\frac34r_{i_{n-1}}}\\
      &\ll& q_{i_{n-1}}^{-2}.
      \eeas
      It follows that $T^{j_{i_n}-q_{i_{n-1}}}x\in I_{i_{n-1}}$ and hence $j_{i_n}-q_{i_{n-1}}\geq j_{i_{n-1}}$.
      This contradicts with $j_{i_n}-j_{i_{n-1}}< q_{i_{n-1}}$.

      Similarly, if $j_{-i_n}-j_{-i_{n-1}}< q_{i_{n-1}}$, then we have $j_{-i_n}\leq q_{i_{n-1}}$. Then it can be reduced to the case with $x\in I_{i_{n}}$.  By Lemma \ref{lemma.specialcase}, we obtain the estimates (\ref{keylemma.norm}) and (\ref{keylemma.angle}).

  \item It is left to consider the case with the following:
  %$$j_{i_n},\ \  j_{-i_n}\geq 10^{-3}q_{i_n},$$
  \beas
  &j_{i_n}-j_{i_{n-1}},\ \  j_{-i_n}-j_{-i_{n-1}}\geq q_{i_{n-1}},\\
  &j_{i_{n-1}},\ \ j_{-i_{n-1}}\geq 10^{-3}q_{i_n}.
  \eeas
  See Figure 4.

      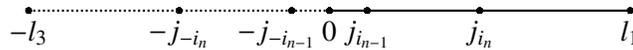
\begin{figure}[H]
\begin{center}
\begin{tikzpicture}[yscale=1]
\draw [densely dotted, thick] (-4,0) -- (0,0);
\draw [thick] (0,0) -- (4,0);
\draw plot [mark=*, mark size= 0.2ex] (-4,0);
\node[align=left, below] at (-4,0){$-l_3$};
\draw plot [mark=*, mark size= 0.2ex] (0,0);
\node[align=left, below] at (0,0){$0$};
\draw plot [mark=*, mark size= 0.2ex] (4,0);
\node[align=left, below] at (4,0){$l_1$};
\draw plot [mark=*, mark size= 0.2ex] (0.5,0);
\node[align=left, below] at (0.5,0){$j_{i_{n-1}}$};
\draw plot [mark=*, mark size= 0.2ex] (2,0);
\node[align=left, below] at (2,0){$j_{i_{n}}$};
\draw plot [mark=*, mark size= 0.2ex] (-0.5,0);
\node[align=left, below] at (-0.7,0){$-j_{-i_{n-1}}$};
\draw plot [mark=*, mark size= 0.2ex] (-2,0);
\node[align=left, below] at (-2,0){$-j_{-i_{n}}$};
\end{tikzpicture}
\caption{The positions of the return times}
\end{center}
\end{figure}

First, we consider $A_{j_{i_n}}(x)=A_{j_{i_{n}}-j_{i_{n-1}}}(T^{j_{i_{n-1}}}x)A_{j_{i_{n-1}}}(x)$.
Since $j_{i_{n-1}}\geq 10^{-3}q_{i_n}\gg q_{i_{n-1}}$ and $j_{-i_{n-1}}\geq 10^{-3}q_{i_n}\gg q_{i_{n-1}}$, we can further assume that $j_{i_{n-1}}\geq r_{i_{n-1}}^+$ and $j_{-i_{n-1}}\geq r_{i_{n-1}}^-$.
By Lemma \ref{claim.1}, it holds that
\bea&\|A_{j_{i_{n-1}}}(x)\|\geq \lambda^{(1-10^{-4})j_{i_{n-1}}},\quad\|A_{-j_{-i_{n-1}}}(x)\|\geq \lambda^{(1-10^{-4})j_{-i_{n-1}}},\nn\\
\label{keylemma.angle.1}
&|s_{j_{i_{n-1}}}(x)-s_{r^+_{i_{n-2}}}(x)| <\lambda^{-\frac32 r_{i_{n-2}}},\\
&|u_{j_{-i_{n-1}}}(x)-u_{r^-_{i_{n-2}}}(x)| <\lambda^{-\frac32 r_{i_{n-2}}}.\label{keylemma.angle.4}
\eea
Since $x\in I_{i_{n-2}}\backslash I_{i_{n-1}}$, then we have $\|x-C^{(i_{n-1})}\|\geq q_{i_{n-1}}^{-2}$. Combining (\ref{keylemma.angle.1}) and (\ref{keylemma.angle.4}), it holds that
\beas
|s_{j_{i_{n-1}}}(x)-u_{j_{-i_{n-1}}}(x)|
\geq |g_{i_{n-2}+1}(x)|-\lambda^{-\frac32 r_{i_{n-2}}}
\geq C q_{i_{n-1}}^{-6}.
\eeas
Then one can apply Lemma \ref{lemma.basic} to obtain
\beq\label{keylemma.angle.5}
|u_{j_{i_{n-1}}}(T^{j_{i_{n-1}}}x)-u_{r^-_{i_{n-1}}}(T^{j_{i_{n-1}}}x)| < \|A_{j_{i_{n-1}}}(x)\|^{-\frac32}\leq \lambda^{- 10^{-6}q_{i_n}}.
\eeq
If $j_{i_n}-j_{i_{n-1}}<10^{-6}q_{i_n}$, then $$\|A_{j_{i_n}-j_{i_{n-1}}}(T^{j_{i_{n-1}}}x)\|\leq \lambda^{10^{-6}q_{i_n}}\ll \lambda^{(1-10^{-4})j_{i_{n-1}}}\leq \|A_{j_{i_{n-1}}}(x)\|.$$
Hence by Lemma \ref{lemma.basic}, it holds that
\bea\label{keylemma.angle.6_1}
\| A_{j_{i_n}}(x)\| &\geq & \|A_{j_{i_{n-1}}}(x)\|\cdot\|A_{j_{i_n}-j_{i_{n-1}}}(T^{j_{i_{n-1}}}x)\|^{-1}\\ &\geq&\lambda^{(1-10^{-4})j_{i_{n-1}}-10^{-6}q_{i_n}}\nn \\&\geq& \lambda^{(1-2\cdot10^{-4})j_{i_n}},\nn
\eea
\bea
\label{keylemma.angle.6}
&|s_{j_{i_{n}}}(x)-s_{j_{i_{n-1}}}(x)| \leq \|A_{j_{i_{n-1}}}(x)\|^{-\frac32}<\lambda^{-10^{-6}q_{i_n}}. \eea
Otherwise, we have $j_{i_n}-j_{i_{n-1}}\geq 10^{-6}q_{i_n}$, which from Lemma \ref{claim.1} implies that
\bea
&\|A_{j_{i_{n-1}}}(x)\|\geq \lambda^{(1-10^{-4})j_{i_{n-1}}},\ \ \|A_{j_{i_{n}}-j_{i_{n-1}}}(T^{j_{i_{n-1}}}x)\|\geq \lambda^{(1-10^{-4})(j_{i_{n}}-j_{i_{n-1}})}.
\nn\\\label{keylemma.angle.7}
&|s_{j_{i_n}-j_{i_{n-1}}}(T^{j_{i_{n-1}}}x) - s_{r^+_{i_{n-1}}}(T^{j_{i_{n-1}}}x)|< \lambda^{-\frac32 r_{i_{n-1}}}.
\eea
Since $T^{j_{i_{n-1}}}x\in I_{i_{n-1}}\backslash I_{i_{n}}$, then we have $\|T^{j_{i_{n-1}}}x-C^{(i_{n})}\|\geq q_{i_{n}}^{-2}$.
Combining (\ref{keylemma.angle.5}) and (\ref{keylemma.angle.7}), we obtain
\beas
&&|u_{j_{i_{n-1}}}(T^{j_{i_{n-1}}}x) - s_{j_{i_n}-j_{i_{n-1}}}(T^{j_{i_{n-1}}}x)| \\&\geq & |g_{i_{n-1}+1}(T^{j_{i_{n-1}}}x)|-\lambda^{-\frac32 r_{i_{n-1}}}-\lambda^{-10^{-6}q_{i_n}}\\ &\geq& Cq_{i_n}^{-6}.
\eeas
Again one can apply Lemma \ref{lemma.basic} to obtain (\ref{keylemma.angle.6_1}) and (\ref{keylemma.angle.6}).

Next we consider $A_{l_1}(x)=A_{j_{i_n}-j_{i_{n}}}(T^{j_{i_{n}}}x)A_{j_{i_{n}}}(x)$.
Similar to the above, we have
\bea
&\| A_{-j_{-i_n}}(x) \| \geq \lambda^{(1-2\cdot 10^{-4})j_{-i_n}},\nn\\
\label{keylemma.angle.8}
&|u_{j_{-i_n}}(x) - u_{j_{-i_{n-1}}}(x)|< \lambda^{- 10^{-6}q_{i_n}}.
\eea
Hence by (\ref{keylemma.angle.1}), (\ref{keylemma.angle.4}), (\ref{keylemma.angle.6}) and (\ref{keylemma.angle.8}), it holds that
\beas
|s_{j_{i_{n}}}(x)-u_{j_{-i_n}}(x)|
&\geq& | s_{j_{i_{n-1}}}(x) - u_{j_{-i_{n-1}}}(x) |- \lambda^{- 10^{-6}q_{i_n}} \\
&\geq& |g_{i_{n-1}+1}(x)|- \lambda^{- 10^{-6}q_{i_n}}\\
&\geq& C q_{i_{n-1}}^{-6}.
\eeas
Again using Lemma \ref{lemma.basic}, it holds that
\beq\label{keylemma.angle.11}
|u_{j_{i_n}}(T^{j_{i_n}}x)-u_{r^-_{i_{n}}}(T^{j_{i_n}}x)| <\|A_{j_{i_n}}(x)\|^{-\frac32}< \lambda^{- 10^{-6}q_{i_n}}.
\eeq
If $l_1-j_{i_n}<10^{-6}q_{i_n}$, then $$\|A_{l_1-j_{i_n}}(T^{j_{i_{n}}}x)\| \leq \lambda^{10^{-6}q_{i_n}} \ll \|A_{j_{i_n}}(x)\|.$$
Hence by Lemma \ref{lemma.basic}, it holds that
\bea
&\| A_{l_1}(x) \| \geq \|A_{j_{i_n}}(x)\|\cdot \|A_{l_1-j_{i_n}}(T^{j_{i_{n}}}x)\|^{-1} \geq\lambda^{(1-2\cdot10^{-4})j_{i_n}-10^{-6}q_{i_n}} \geq \lambda^{\frac{99}{100}l_1},
\nn\\\label{keylemma.angle.9}
&|s_{l_1}(x)-s_{j_{i_n}}(x)|\leq \|A_{j_{i_n}}(x)\|^{-\frac32}<\lambda^{- 10^{-6}q_{i_n}},
\eea
where the first estimate is (\ref{keylemma.norm}).
Otherwise, we have $l_1-j_{i_n}\geq 10^{-6}q_{i_n}$, which by Lemma \ref{claim.1} implies that
\beq\label{keylemma.angle.13}
|s_{l_1-j_{i_n}}(T^{j_{i_n}}x) - s_{r^+_{i_n}}(T^{j_{i_n}}x)| <\lambda^{-\frac32r_{i_n}}.
\eeq
Since $T^{j_{i_n}}x\in I_{i_n}\backslash I_{i_{n+1}}$, then $\|T^{j_{i_n}}x-C^{(i_{n+1})}\|\geq q_{i_{n+1}}^{-2}$. Then using (\ref{keylemma.angle.11}) and (\ref{keylemma.angle.13}), we obtain
\beas
&&|u_{j_{i_n}}(T^{j_{i_n}}x)-s_{l_1-j_{i_n}}(T^{j_{i_n}}x)|\\
&\geq& |g_{i_n+1}(T^{j_{i_n}}x)|-\lambda^{-\frac32r_{i_n}}-\lambda^{- 10^{-6}q_{i_n}}\\
&\geq& Cq_{i_{n+1}}^{-6}.
\eeas
By (\ref{keylemma.angle.6_1}) and Lemma \ref{claim.1}, we have
\beqs
\|A_{j_{i_n}}(x)\|\geq \lambda^{(1-2\cdot10^{-4})j_{i_n}},\ \ \|A_{l_1-j_{i_n}}(x)\|\geq \lambda^{(1-10^{-4})(l_1-j_{i_n})}.
\eeqs
Again one can apply Lemma \ref{lemma.basic} to obtain (\ref{keylemma.norm}) and (\ref{keylemma.angle.9}). Similarly, we obtain
\bea\nn
&\| A_{-l_3}(x)\| \geq \lambda^{\frac{99}{100}l_3},
\\\label{keylemma.angle.10}
&|u_{l_3}(x)-u_{j_{-i_{n}}}(x)| <\lambda^{-10^{-6}q_{i_n}}.
\eea
Combining (\ref{keylemma.angle.6}), (\ref{keylemma.angle.8}), (\ref{keylemma.angle.9}) and (\ref{keylemma.angle.10}), we have that
\bea\label{keylemma.angle.12}
|s_{l_1}(x)-u_{l_3}(x)| &\geq&  |s_{j_{i_n}}(x)-u_{j_{-i_n}}(x)|-\lambda^{-10^{-6}q_{i_n}}\\& \geq& |s_{j_{i_{n-1}}}(x)-u_{j_{-i_{n-1}}}(x)|-2\lambda^{-10^{-6}q_{i_n}}\nn\\
&\geq& |g_{i_{n-2}+1}(x)|-\lambda^{-\frac32 r_{i_{n-2}}}-2\lambda^{-10^{-6}q_{i_n}}\nn\\
&\geq& C q_{i_{n-1}}^{-6}.\nn
\eea
Hence by Lemma \ref{lemma.basic} the estimate (\ref{keylemma.angle}) follows.

\end{enumerate}

At the end, we deal with the cases with $q_{i_{n}+1}\leq q_{i_{n}}^3$ or $q_{i_{n-1}+1}\leq q_{i_{n-1}}^3$. It is easy to see from Corollary \ref{claim.3} that they are degenerate situations in the case with $q_{i_{n}+1}> q_{i_{n}}^3$ and $q_{i_{n-1}+1}> q_{i_{n-1}}^3$. The proofs are nearly the same and we will not repeat them here.

Until now we have completed the proof of Lemma \ref{keylemma}. \hfill$\Box$\vspace{0.7mm}

\

\noindent
\begin{prof2}
If $T^ix\in I_{i_{n}}$, then (\ref{theorem.LDT.5}) follows from the same argument for (\ref{lemma.specialcase.0});
if $T^ix\in I_{i_{n-1}}\backslash I_{i_{n}}$, then (\ref{theorem.LDT.5}) can be obtained with the same method for (\ref{lemma.specialcase.6});
if $T^ix\in I_{i_{n-2}}\backslash I_{i_{n-1}}$, then (\ref{theorem.LDT.5}) is similar to the proof of (\ref{keylemma.angle.12}).
\end{prof2}

\section{Proof of Theorem \ref{Holder}}
%We will give the proof of the H\"older continuity of the LE and the IDS in this section.
We first introduce the Avalanche Principle, which is given by Goldstein-Schlag \cite{GS1}.
\begin{theorem}\label{AP}\textbf{(Avalanche Principle)}
Let $B_1,\cdots,B_m$ be a sequence in $\text{SL}(2,\mathbb{R})$ satisfying
\beq\label{8} \min_{1\leq j \leq m} \|B_j \|\geq \mu >m,\eeq
\beq\label{9} \max_{1\leq j \leq m} \left|\log \|B_j \|+\log \|B_{j+1} \|-\log \|B_{j+1}B_{j} \| \right|<\frac12 \log \mu.\eeq
Then it holds that
\beq\label{10}
\left|\log \|B_m \cdots B_1\|+\sum_{j=2}^{m-1}\log\|B_j\| -\sum_{j=1}^{m-1}\log\|B_{j+1}B_j\|\right|<C\frac m \mu.
\eeq
\end{theorem}

Theorem \ref{Holder} follows from Theorem \ref{LDT} and the Avalanche Principle by some inductive steps in \cite{GS1}.
We are going to rewrite the proof here because we want to show the H\"older exponent $\sigma$ can be chosen to be independent of $\lambda$ from (\ref{LDT.LDT}). This improves the results in \cite{GS1,YZ}, where the H\"older exponent gets worse as $\lambda$ gets larger.

\begin{prof1}
Choose $n_1$ large enough such that Theorem \ref{LDT} holds.
Let $n_2=mn_1=\lambda^{10^{-7}n_1}n_1$. Then Theorem \ref{LDT} permits us to ensure that
\beq\label{6}
n_1\log\lambda\geq \log\|A_{n_1}(T^{jn_1}x)\| \geq \frac56n_1\log\lambda ,\eeq
\beq
\label{7}  2n_1\log\lambda\geq \log\|A_{2n_1}(T^{2jn_1}x)\| \geq\frac{10}6n_1\log\lambda ,\eeq
for $j=0,1,2,\cdots,m$, except for $x\in \Omega \subseteq \R/\Z$, where
$$\text{meas}\ \Omega <4n_2\lambda^{-10^{-6}n_1}<\lambda^{-10^{-7}n_1}.$$
Fix $x\notin \Omega$, and define
$$B_j=A_{n_1}(T^{jn_1}x),\quad j=0,1,\cdots,m.$$
We verify the conditions of Lemma \ref{AP}. Let $\mu= \lambda^{\frac56n_1}$.
By (\ref{6}) and (\ref{7}),
$$\|B_j\|\geq \lambda^{\frac56n_1}=\mu>m,$$
$$\left| \log\|B_j\|+\log\|B_{j+1}\|-\log\|B_{j+1}B_{j}\| \right|
\leq \frac13n_1\log\lambda<\frac12\log\mu.$$
Thus (\ref{8}) and (\ref{9}) hold, and (\ref{10}) implies
$$ \left|\log \|A_{n_2}(x) \| +\sum_{j=2}^{m-1}\log \|A_{n_1}(T^{jn_1}x) \| -\sum_{j=1}^{m-1}\log \|A_{2n_1}(T^{jn_1}x) \| \right|<C\frac m \mu.$$
This leads to
$$\left|L_{n_2}(E)+\frac{m-2}{m}L_{n_1}(E)-2\frac{m-1}mL_{2n_1}(E)\right|<\frac C{n_1\mu}+ \log\lambda \cdot \text{meas} \  \Omega <\lambda^{-10^{-7}n_1}.$$
Then
$$\left|L_{n_2}(E)+L_{n_1}(E)-2L_{2n_1}(E)\right|<C\lambda^{-10^{-7}n_1}.$$
The same inequality also holds true with $n_2$ replaced by $2n_2$. Thus we obtain
$$\left|L_{2n_2}(E)-L_{n_2}(E)\right|<C\lambda^{-10^{-7}n_1}.$$
Iterating the above procedures, we have
$$\left|L_{n_{s+1}}(E)+L_{n_s}(E)-2L_{2n_s}(E)\right| <C\lambda^{-10^{-7}n_s},$$
$$\left|L_{2n_{s+1}}(E)-L_{n_{s+1}}(E)\right|<C\lambda^{-10^{-7}n_{s}},  \ \  s\geq 1,$$
where $n_{s+1}=\lambda^{10^{-7}n_s}n_s$. Then
\beas
\left|L(E)-L_{n_{2}}(E)\right|&\leq& \sum_{s\geq 2}\left|L_{n_{s+1}}(E)-L_{n_s}(E)\right|\\
&\leq& \sum _{s\geq 2}\left(2\left|L_{2n_{s}}(E)-L_{n_s}(E)\right|+C\lambda^{-10^{-7}n_s}\right)\\
&\leq& \sum _{s\geq 2}C\lambda^{-10^{-7}n_s} <C\lambda^{-10^{-7}n_2}.
\eeas
And hence
\bea&&
\left|L(E)+L_{n_{1}}(E)-L_{2n_{1}}(E)\right|\label{L}\\&\leq& \left|L(E)-L_{n_{2}}(E)\right| \nn +\left|L_{n_{2}}(E)+L_{n_{1}}(E)-L_{2n_{1}}(E)\right| \\ &\leq& C\lambda^{-10^{-7}n_1}.\nn
\eea
Obviously, for any $n$,
$$\left\|\partial_{E}A_n(x)\right\|\leq n \left(C\lambda+2\right)^{n-1}.$$
And then
\beq\label{Ln1}\left|L_{n_{1}}(E)-L_{n_{1}}(E')\right|\leq \left(C\lambda\right)^{n_{1}}\left|E-E'\right|,\eeq
\beq\label{Ln2}\left|L_{2n_{1}}(E)-L_{2n_{1}}(E')\right|\leq \left(C\lambda\right)^{2n_{1}}\left|E-E'\right|.\eeq
Combining (\ref{L}), (\ref{Ln1}) and (\ref{Ln2}),
$$\left|L(E)-L(E')\right| <\left(C\lambda\right)^{2n_{1}}\left|E-E'\right| +C\lambda^{-10^{-7}n_1}.$$
Since the choice of $n_1$ is arbitrary as long as it is large enough, then we let
$$n_1= -\frac{\log|E-E'|}{4\log \lambda}.$$
Therefore
\beqs  \left|L(E)-L(E')\right|<C\left|E-E'\right|^{\sigma},\eeqs
where $\sigma=10^{-8}$.
Using the Thouless' Formula and the Hilbert transform, we can also obtain the H\"older continuity of the IDS.
This proves the theorem.
\end{prof1}

\section{Appendix}
\subsection{Arithmetical properties of the frequencies}
Let $I_n=B\left(c,q_{n}^{-2}\right)$, where $c$ is a constant in $\R/\Z$ and $\left\{{p_n}/{q_n}\right\}$ is the continued fraction approximation of $\alpha$.
It is well known that $\left\{{p_n}/{q_n}\right\}$ satisfies
\beq\label{arithmetical.1}\frac1{q_n(q_{n+1}+q_n)} \leq\left|\alpha-\frac{p_n}{q_n}\right|\leq\frac1{q_n q_{n+1}},\eeq
\beq\label{arithmetical.2}\min_{0<k<q_n}\|k\alpha\| \geq \|q_{n-1}\alpha\| \geq \frac1{q_{n-1}+q_{n}} .\eeq
By (\ref{arithmetical.1}), it holds that
$$\left|k\alpha-k\frac{p_n}{q_n}\right|\leq\frac k{q_nq_{n+1}}\leq\frac1{q_{n+1}},$$
for $0\leq k < q_n$.
This means that $\alpha$ is extremely close to the rational number ${p_n}/{q_n}$.
Since the distribution of $\left\{\left\|\frac{k p_n}{q_n}\right\|:0\leq k< q_n\right\}$ in [0,1] is orderly arranged with a step size $q_n^{-1}$, then so is $\left\{k\alpha : 0\leq k< q_n\right\}$ with an error of $q_{n+1}^{-1}$.

The next lemma give the basic property of the return time of the trajectory.
\begin{lemma}\label{lemma.returntime}
For $x\in I_n$, the return time of the trajectory to $I_n$ is not smaller than $q_n$.
\end{lemma}
\begin{prof}
By (\ref{arithmetical.2}), we have that
\beas
\min_{0<k<q_n}\|k\alpha\|\geq {(2q_n)}^{-1}\geq 2q_n^{-2}.
\eeas
If $x\in I_n$, then $x+k\alpha\notin I_n$ for $0<|k|<q_n$.
\end{prof}

\

\subsection{Types of the angle functions}
In Theorem \ref{theorem.induction}, the angle function possesses different properties in different cases.
The angle function $g_{i+1}$ behaves as an affine function in Case $(i+1)_\I$, as an quadratic function in Case $(i+1)_\II$, and as (\ref{gi+1-gi3}) in Case $(i+1)_\III$, which is called type $\I$ function, type $\II$ function and type $\III$ function, respectively.
In this subsection, we would like to give the definition and some properties. One can refer to \cite{WZ} for the proofs of these lemmas.

\begin{definition}\label{def.types}
For a connected interval $J\subset \R/\Z$ and constant $0<a\leq 1$, let $aJ$ be the subinterval of $J$ with the same center and whose length is $a|J|$. Let $I_n$ and $g_{n+1}$ be as in Theorem \ref{theorem.induction}. Let $I=I_{n,1}$ or $I_{n,2}$. We say
\begin{itemize}
  \item  \textbf{$g_{n+1}$ is of type $\I$ in $I$}, if %$I_{n,1}\cap I_{n,2} = \varnothing$ and for $j=1,2$ ,
          \begin{itemize}
            \item $\|g_{n+1}(x)\|_{C^2}<C\lambda^{2r_{n-1}}$ and $g_{n+1}(x)=0$ has only one solution in $I$, say $c_{n+1}$, which is contained in $\frac13I$;
            \item $\frac {dg_{n+1}}{dx}(x)=0$ has at most one solution in $I$ while $\left|\frac{dg_{n+1}}{dx}\right|>q_n^{-4}$ for all $x\in \frac 1 2 I$;
            \item Let $J\subset I$ be the subinterval such that $\frac{dg_{n+1}}{dx}(J)\frac{dg_{n+1}}{dx}(c_{n+1})\leq 0$, then $|g_{n+1}(x)|>cq_n^{-6}$ for all $x\in J$.
          \end{itemize}
          Let $\I_+$ denotes the case $\frac {dg_{n+1}}{dx}(c_{n+1,j})>0$ and $\I_-$ for $\frac {dg_{n+1}}{dx}(c_{n+1,j})<0$. We also say $g_n$ is of type $\I$ in $I_{n,1}$ or $I_{n,2}$, if we replace $g_{n+1}$ with $g_n$ in the above definition.
  \item \textbf{$g_{n+1}$ is of type ${\II}$ in $I$}, if
           \begin{itemize}
             \item $\|g_{n+1}\|_{C^2}<C$ and $g_{n+1}(x)=0$ has at most two solutions, say $c_{n+1,1}$ and $c_{n+1,2}$, which are contained in $\frac12I$;
             \item $\frac {dg_{n+1}}{dx}(x)=0$ has one solution which is contained in $\frac12I$;
             \item $ g_{n+1}(x)=0$ has only one solution if and only if  $\frac {dg_{n+1}}{dx}(c_{n+1,1})=0$;
             \item $\left|\frac {d^2g_{n+1}}{dx^2}\right|>c$ whenever $\left|\frac {dg_{n+1}}{dx}\right|<q_n^{-4}$.
           \end{itemize}
  \item \textbf{$g_{n+1}$ is of type ${\III}$ in $I$}, if
           \beas
           &g_{n+1}(x)=\arctan(l_k^2 \tan(f_1(x)))-\frac\pi2+f_2(x),\\
           &\|f_1(x)-g_{n}(x\pm k\alpha)\|_{C^2}<C\lambda^{-\frac32r_{n-1}},\\
           &\|f_2(x)-g_{n}(x)\|_{C^2}<C\lambda^{-\frac32r_{n-1}},
           \eeas
           where $l_k=\|A_{k}(x)\|$ with $|k|<10^{-6}q_n$, and "$\pm$" takes "$+$" for $x\in I_{n,1}$ and takes "$-$" for $x\in I_{n,2}$. Moreover, either $g_n(\cdot\pm k\alpha)$ is of type $\I_+$ and $g_n(\cdot)$ is of type $\I_-$, or $g_n(\cdot\pm k\alpha)$ is of type $\I_-$ and $g_n(\cdot)$ is of type $\I_+$.
\end{itemize}
\end{definition}
\begin{figure}[H]
\begin{center}
\begin{tikzpicture}[yscale=1.5]
\draw [->] (-1.5,0) -- (1.5,0);
\draw [->] (0,-0.5) -- (0,1.3);
\draw [thick,domain=-1.4:1.4] plot (\x, {0.1*\x*\x+0.3*\x});
\node [above] at (0,1.3) {$f(x)$};
\node [right] at (1.5,0) {$x$};
\node[align=left, below] at (0,-0.5)%
{Type $\I$};
\end{tikzpicture}\ \ \ \
\begin{tikzpicture}[yscale=1.5]
\draw [->] (-1.5,0) -- (1.5,0);
\draw [->] (0,-0.5) -- (0,1.3);
\draw [thick,domain=-1.4:1.4] plot (\x, {0.25*\x*\x-0.1});
\node [above] at (0,1.3) {$f(x)$};
\node [right] at (1.5,0) {$x$};
\node[align=left, below] at (0,-0.5)%
{Type ${\II}$};
\end{tikzpicture}\ \ \ \
\begin{tikzpicture}[yscale=0.7]
\draw [->] (-1.5,0.6) -- (1.5,0.6);
\draw [->] (0,0) -- (0,3.75);
\draw [thick,domain=-1:1] plot (\x, {1/180*pi*atan(20*tan(\x r))+0.5*pi-0.5*(\x-0.9)});
\draw [semithick,domain=-1.2:1.2] plot (\x, {0*\x+pi});
\node [above] at (0,3.75) {$f(x)$};
\node [right] at (1.5,0.6) {$x$};
\node [left] at (-1.5,0.6) {$-\pi$};
\node [left] at (-1.2,3.14) {$0$};
\node[align=left, below] at (0,-0.1)%
{Type ${\III}$};
\end{tikzpicture}
\caption{Graphs of Type $\I$, ${\II}$ and ${\III}$ functions}
\label{graph}
\end{center}
\end{figure}
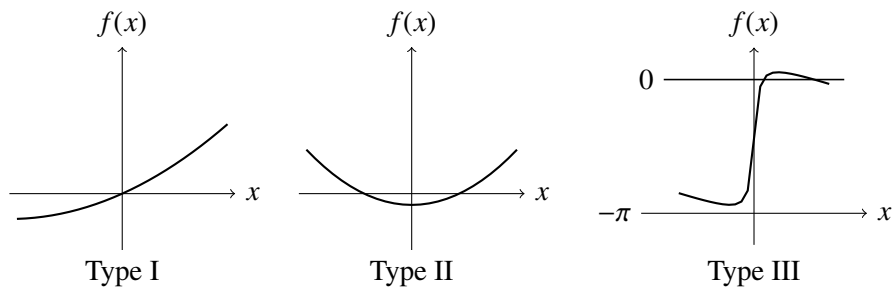
The functions of type $\I$ behaves as an affine function and the functions of type $\II$ behaves as a quadratic function. See Figure 5.
For the functions of type ${\III}$, they can be divided into two parts: $\arctan \left(l_k^2\tan g_n(x\pm k\alpha)\right)-\frac\pi 2$ and $g_n(x)$.
If $x$ locates far away from the zero of $g_n(\cdot\pm k\alpha)$, then the first part is so small to be neglected and therefore the function $f$ approximately equals to $g_n(\cdot)$. In this case $g_{n+1}$ is in fact of type $\I$.
If $x$ varies near the zero of $g_n(\cdot\pm k\alpha)$, the first part has a drastic change from $-\pi $ to $0$, which leads to a bifurcation of no zero, one zero and two zeros of $g_{n+1}$.
The next lemma gives a careful description of the functions of type $\III$.
\begin{lemma}\label{lemma.propertyoftpye3}
Let $g_{n+1}$ be of type $\III$. Without loss of generality, $g_n(\cdot\pm k\alpha)$ is of type $\I_+$ and $g_n(\cdot)$ is of type $\I_-$. Assume  $f_1(e_1)=f_2(e_2)=0$ and $d:=|c_{n,1}+k\alpha-c_{n,2}|\leq \frac 23q_{n}^{-2}$. Consider $x\in I_{n,1}$. Let
$$X=\left\{x\in I_{n,1}:\mathbb{RP}^1\ni\left|g_{n+1}(x)\right|=\min_{y\in I_{n,1}}|g_{n+1}(y)|\right\}=\{c_{n+1,1},c'_{n+1,2}\},$$
with $|c_{n+1,1}-c_{n,1}|<C\lambda^{-\frac34 r_{n-1}}$.
Let $\eta_j$ be constants satisfying $q_n^{-2}\leq \eta_j\leq q_n^{2}$, $0\leq j\leq 4$. Then
\beq\label{lemma.typeIII.1}\left|e_{1}-c'_{n+1,2}\right|<Cl_k^{-\frac34},\ \ |e_2-c_{n+1,1}|<Cl_k^{-\frac34}.\eeq
\iffalse In particular, if $g_{n+1}(X)=\{0\}$, then
$$0 < c'_{n+1,2} \leq x_2 < d;$$
if $f(x_1)=f(x_2)\neq 0$, then
$$x_1=x_2.$$\fi
Moreover, there exist two distinct  points $x_1,x_2\in B\left(c'_{n+1,2}, \eta_0 l_k^{-1}\right)$ such that
$$\frac{dg_{n+1}}{dx}(x_j)=0,\quad j=1,2,$$
and $x_1$ is a local minimum with $$g_{n+1}(x_1)>\eta_1l_k^{-1}-\pi.$$
Furthermore, let $l_{k}^{-\frac14}\leq r\leq\frac14q_{n}^{-2}$, we have the following.\\
If $d\geq \frac 13 q_n^{-2}$, then we have
$$\left|g_{n+1}(x)\right|>cr^3,\quad x\notin B\left(c'_{n+1,2},Cl_k^{-\frac14}\right)\cup B\left(c_{n+1,1},r\right);$$
$$ \left\|g_{n+1}-f_2\right\|_{C^2}<Cl_k^{-\frac32},\quad x\in B\left(c_{n+1,1},\frac14q_n^{-2}\right).$$
If $d\leq \frac 13 q_n^{-2}$, then we have
$$|g_{n+1}(x)|>cr^3,\quad x\notin B\left(X, \frac r6\right);$$
$$\left|\frac {d^2g_{n+1}}{dx^2}(x)\right|>c \ \ \text{whenever}\ \ \left|\frac {dg_{n+1}}{dx}(x)\right|<q_n^{-4},\ \  x\in B\left(X,\frac 16q_{n}^{-2}\right).$$
Finally, we have the following bifurcation as d varies. There is a $d_0=\eta_2 l_k^{-1}$ such that
\begin{itemize}
  \item if $d>d_0$, then $g_{n+1}(x)=0$ has two solutions.
  \item If $d=d_0$, then $g_{n+1}(x)=0$ has exactly one tangential solution. In other words, $c_{n+1,1}=c'_{n+1,2}=x_2$ and $g_{n+1}(x_2)=0$.
  \item If $d<d_0$, then $g_{n+1}(x)\neq 0$ for $x\in I_{n,1}$. Moreover, $\min_{x\in I}|g_{n+1}(x)| =-\eta_3 l_k^{-1}+\eta_4 d$.
\end{itemize}
\end{lemma}

\begin{prof}
Let $f=g_{n+1}$ in Lemma 6 in \cite{WZ}.
\end{prof}

\begin{lemma}\label{typeIII}
If the $(i+1)$th step are in Case $(i+1)_{\I}$, then in the $i$th step the angle function $g_{i}$ must be of type $\I$  in $I_{i}$;
If the $(i+1)$th step are in Case $(i+1)_{\II}$, then in the $i$th step the angle function $g_{i}$ must be of type $\II$  in $I_{i-1}$, or of type $\III$ in $I_{i-1}$;
if the $(i+1)$th step are in Case $(i+1)_{\III}$, then in the $i$th step $g_{i}$ must be of type $\I$ in $I_{i}$, or of type $\III$ in $I_{i-1}$ with
\beqs\|g_{i+1}-g_i\|_{C^2}\leq C\lambda^{-\frac32r_{i-1}}\quad \text{in} \ \ I_i.\eeqs
\end{lemma}
\begin{prof}
It is obvious by the proof of Theorem 3 in \cite{WZ}.
\end{prof}

Note that if $d$ in Lemma \ref{lemma.propertyoftpye3} is sufficiently close to $d_0$ and $r$ is sufficiently small, then the restriction of type $\III$ function $f$ to $B(x_1,r)$ may become type $\II$. In Lemma \ref{typeIII} the process from Case $(i)_\III$ to Case $(i+1)_{\II}$ is actually this restriction.  For simplicity, we still call $f$ is of type $\III$ because all the estimates of type $\III$ are applicable. Hence we do not need to consider type $\II$ except that the angle function is of type $\II$ for every step, which is a trivial case.

\

\noindent\bf{\footnotesize Acknowledgements}\quad\rm
{\footnotesize  The authors are grateful to Zhenghe Zhang and Kristian Bjerkl\"ov for useful suggestions. }

\

 \end{document}